\documentclass[amstex,12pt, amssymb]{article}

\usepackage{mathtext}
\usepackage[cp1251]{inputenc}
\usepackage[T2A]{fontenc}
\usepackage[dvips]{graphicx}
\usepackage{amsmath}
\usepackage{amssymb}
\usepackage{amsxtra}
\usepackage{latexsym}
\usepackage{ifthen}

\def\Xint#1{\mathchoice
   {\XXint\displaystyle\textstyle{#1}}
   {\XXint\textstyle\scriptstyle{#1}}
   {\XXint\scriptstyle\scriptscriptstyle{#1}}
   {\XXint\scriptscriptstyle\scriptscriptstyle{#1}}
   \!\int}
\def\XXint#1#2#3{{\setbox0=\hbox{$#1{#2#3}{\int}$}
     \vcenter{\hbox{$#2#3$}}\kern-.5\wd0}}
\def\dashint{\Xint-}

\newcounter{lemma}[section]

\newcounter{corollary}[section]

\newcounter{remark}[section]

\newcounter{theorem}[section]

\newcounter{proposition}[section]

\newcounter{example}

\numberwithin{equation}{section}

\begin{document}

\markboth{\centerline{V.~GUTLYANSKII, V.~RYAZANOV, E.~SEVOST'YANOV
and E.~YAKUBOV}}{\centerline{ON THE DIRICHLET PROBLEM ...}}

\def\cc{\setcounter{equation}{0}
\setcounter{figure}{0}\setcounter{table}{0}}

\overfullrule=0pt

%\normalsize\large

\author{Gutlyanskii V., Ryazanov V., Sevos'yanov E., Yakubov E.}

\title{Toward the theory of Dirichlet problem \\ for the degenerate Beltrami equations}

\date{\today}
\maketitle

%\large
\begin{abstract}
In this article, first we give a general lemma on the existence of
regular ho\-meo\-mor\-phic solutions $f$ with the hydrodynamic
normalization $f(z)=z+o(1)$ as $z\to\infty$ to the
de\-ge\-ne\-ra\-te Beltrami equations
$\overline{\partial}f=\mu\,\partial f$ in $\mathbb C$ whose complex
coefficients $\mu$ have compact supports. On this basis, we
establish criteria for existence and representation of regular
discrete open solutions for the Dirichlet problem with continuous
data to degenerate Beltrami equations in arbitrary simply connected
bounded domains $D$ in $\mathbb C$.

%\medskip

Moreover, we obtain similar criteria for the existence of
multi-valued solutions $f$ in the spirit of the theory of
multi-valued analytic functions in arbitrary bounded domains $D$ in
$\mathbb C$  with no boundary component degenerated to a single
point. Note that the latter request is necessary and that the real
parts $u$ of such solutions $f$ are the so-called $A-$harmonic
functions, i.e., single-valued continuous weak solutions of elliptic
equations ${\rm div}\, (A\nabla u)=0$ with matrix-valued
coefficients $A$ associated with $\mu$. Thus, the results can be
applied to potential theory in anisotropic and inhomogeneous media.

%\medskip

The corresponding criteria for the existence of the homeomorphic
solutions and solutions of the Dirichlet problems are formulated in
terms of the quantity
$$
K_{\mu}(z)\ =\ \frac{1+|\mu(z)|}{1-|\mu(z)|}\ ,\ \ \ \ \ \ z\in D\ ,
$$
as well as in terms of the more refined quantities
$$
K^T_{\mu}(z,z_0)\ =\ \frac{\left|1-\frac{\overline{z-z_0}}{z-z_0}\mu
(z)\right|^2}{1-|\mu (z)|^2}\ ,\ \ \ \ \ \ z\in D\ ,\ z_0\in\mathbb
C\ ,
$$
that takes into account not only the modulus of $\mu$ but also its
argument.
\end{abstract}

\bigskip


\bigskip


\bigskip

{\bf Dedicated to the memory of Professor Uri Srebro ( 1936 - 2016
)}

%\newpage

\section{Introduction}

Let $D$ be a domain in the complex plane ${\Bbb C}$, i.e., a
connected open subset of ${\Bbb C}$, and let $\mu:D\to{\Bbb C}$ be a
measurable function with $|\mu(z)|<1$ a.e. (almost everywhere) in
$D$. A {\bf Beltrami equation} is an equation of the form
\begin{equation}\label{eqBeltrami} f_{\bar z}=\mu(z)\,f_z\end{equation}
with the formal complex derivatives $f_{\bar
z}=\overline{\partial}f=(f_x+if_y)/2$, $f_{z}=\partial
f=(f_x-if_y)/2$, $z=x+iy$, where $f_x$ and $f_y$ are partial
derivatives of $f$ in $x$ and $y$, correspondingly. The function
$\mu$ is  said to be the {\bf complex coefficient} and
\begin{equation}\label{eqKPRS1.1}K_{\mu}(z)\ :=\ \frac{1+|\mu(z)|}{1-|\mu(z)|}\end{equation}
the {\bf dilatation quotient} of the equation (\ref{eqBeltrami}).
The Beltrami equation is called {\bf degenerate} if ${\rm ess}\,{\rm
sup}\,K_{\mu}(z)=\infty$.

It is well known that if $K_{\mu}$ is bounded, then the Beltrami
equation has homeomorphic solutions, see e.g. the monographs
\cite{Ahl}, \cite{LV} and \cite{Vek}. Recently, a series of
effective criteria for the existence of homeomorphic solutions have
been also established for degenerate Beltrami equations, see e.g.
historic comments with relevant references in monographs the
\cite{AIM}, \cite{GRSY} and \cite{MRSY}, in BMO-article \cite{RSY}
and in the surveys \cite{GRSY$_*$} and \cite{SY}.

These criteria were formulated both in terms of $K_{\mu}$ and the
more refined quantity that takes into account not only the modulus
of the complex coefficient $\mu$ but also its argument
\begin{equation}\label{eqTangent} K^T_{\mu}(z,z_0)\ :=\
\frac{\left|1-\frac{\overline{z-z_0}}{z-z_0}\mu (z)\right|^2}{1-|\mu
(z)|^2} \end{equation} that is called the {\bf tangent dilatation
quotient} of the Beltrami equation with respect to a point
$z_0\in\mathbb C$, see e.g. \cite{And}, \cite{BGR$_1$},
\cite{BGR$_2$}, \cite{GMSV}, \cite{Le} and
\cite{RSY}-\cite{RSY$_5$}. Note that
\begin{equation}\label{eqConnect} K^{-1}_{\mu}(z)\leqslant K^T_{\mu}(z,z_0) \leqslant K_{\mu}(z)
\ \ \ \ \ \ \ \forall\ z\in D\, ,\ z_0\in \Bbb C\ .\end{equation}
The geometrical sense of $K^T_{\mu}$ can be found e.g. in the
monographs \cite{GRSY} and \cite{MRSY}.

Boundary-value problems  for the Beltrami equations are due to the
famous dissertation of  Riemann (1851) who considered a particular
case of ana\-ly\-ti\-cal  functions when $\mu(z)\equiv 0$, and to
the works of Hilbert (1904, 1912, 1924) who studied the
corres\-pon\-ding system of Cauchy--Riemann for the real and
imaginary parts of analytical functions $f=u+iv$, as well as to the
work of  Poincare (1910) on rising tides.

In this connection, recall that if an analytic function  $f$ defined
in the unit disk ${\Bbb D}= \{z\in {\Bbb C}:$ $ |z| < 1 \}$ is
extended by continuity to its closure, then by the Schwarz formula
\begin{equation}\label{eqDIR4}
f(z)\quad=\quad i\,{\rm Im\,} f(0)\quad+\quad\frac{1}{2\pi i }
\int\limits_{|\zeta|=1} {\rm Re} \,
f(\zeta)\cdot\frac{\zeta+z}{\zeta-z}\frac{d\zeta}{\zeta}\ ,
\end{equation}
see e.g. Section 8, Ch. III, Part  3 in \cite{HC}. Thus, the
analytical  function  $f$ in  ${\Bbb D}$ is determined, up to a
purely imaginary number $ic $, $c={\rm Im\,} f(0),$ by its real part
on the boundary.

Hence the {\bf Dirichlet problem} for the nondegenerate Beltrami
equations in a Jordan domain $D\subset{\Bbb C}$ was the problem on
the existence of a continuous function $f:D\to{\Bbb C}$ with
generalized derivatives by Sobolev of the first order, satisfying
(\ref{eqBeltrami}) a.e., such that
\begin{equation}\label{eqDIRICHLET}\lim\limits_{z\to\zeta}{\rm
Re}\,f(z)=\varphi(\zeta)\qquad\forall\ \zeta\in\partial
D\end{equation} for each prescribed continuous function
$\varphi:\partial D\to{\Bbb R},$ see e.g. \cite{Bojar} and
\cite{Vek}.

Criteria for existence of solutions of the Dirichlet problem in the
unit disk for degenerate Beltrami equations can be found in
monographs \cite{GRSY}, see also survey \cite{GRSY$_*$}. A series of
theorems on the existence of regular, pseudoregular and multi-valued
solutions for the Dirichlet problem to the degenerate Beltrami
equations in Jordan domains and domains bounded by finite
collections of Jordan curves, respectively, was established in
\cite{KPR} and \cite{RSSY}, cf. \cite{GRYY}. Finally, \cite{GRY},
see also \cite{BGR$_3$}, contained the similar results on the
Dirichlet problem to the degenerate Beltrami equations formulated in
terms of prime ends by Caratheodory in arbitrary domains with no
boundary component degenerated to a point.

However, the topology of prime ends has in general a very
complicated nature. Therefore the natural desire arises, in the case
of such general domains, also to obtain the appropriate theorems on
the Dirichlet problem for degenerate Beltrami equations in the usual
sense (\ref{eqDIRICHLET}).

In the present article, we first give a general lemma on the
existence of regular ho\-meo\-mor\-phic solutions $f$ with the
hydrodynamic normalization $f(z)=z+o(1)$ as $z\to\infty$ to the
de\-ge\-ne\-ra\-te Beltrami equations in $\mathbb C$ whose complex
coefficients $\mu$ have compact supports. On this basis, then we
establish criteria for existence and representation of regular
discrete open solutions for the Dirichlet problem with continuous
data to degenerate Beltrami equations in arbitrary simply connected
bounded domains $D$ in $\mathbb C$.

In addition, here we obtain similar criteria for the existence of
multi-valued solutions $f$ in the spirit of the theory of
multi-valued analytic functions in arbitrary bounded domains $D$ in
$\mathbb C$  with no boundary component degenerated to a single
point. Note that the real parts $u$ of such solutions $f$ are the
so-called $A-$harmonic functions, i.e., single-valued continuous
weak solutions of elliptic equations ${\rm div}\, (A\nabla u)=0$
with matrix-valued coefficients $A$ associated with $\mu$. Thus, the
results can be applied to the potential theory in anisotropic and
inhomogeneous media.

Let us emphasize, the request on domains to have no boundary
component degenerated to a single point is necessary. Indeed,
consider the punctured unit disk  $\mathbb D_0:=\mathbb D\setminus\{
0\}$. Setting $\varphi(\zeta)\equiv 1$ on $\partial D$ and
$\varphi(0)= 0$ in (\ref{eqDIRICHLET}), we see that $\varphi$ is
continuous on $\partial\mathbb D_0=\partial\mathbb D\cup\{ 0\}$. Let
us assume that there is an analytic function $f$ satisfying
(\ref{eqDIRICHLET}) with the given $\varphi$. Then the harmonic
function $u:={\rm Re}\, f$ is bounded and by the classic
Cauchy--Riemann theorem, see also Theorem V.4.2 in \cite{Ne}, the
extended $u$ is harmonic in $\mathbb D$. Thus, by contradiction with
the Mean-Value-Property we disprove the above assumption, see e.g.
Theorem 0.2.4 in \cite{ST}.

\section{More definitions and preliminary remarks}

Let $D$ be a domain in the complex plane ${\Bbb C}$. A function
$f:D\to\mathbb C$ in the Sobolev class $W^{1,1}_{\rm loc}$ is called
a {\bf regular solution} of the Beltrami equation (\ref{eqBeltrami})
if $f$ satisfies (\ref{eqBeltrami}) a.e. and its Jacobian $J_f(z)>0$
a.e.

\bigskip

{\bf Lemma 1.} {\it Let a function $\mu : \mathbb C\to{\Bbb C}$ be
with compact support $S$, $|\mu (z)| < 1$  a.e. and $K_{\mu}\in
L^1(S).$ Suppose that, for every $z_0\in S$, there is a family of
measurable functions
${\psi}_{z_0,{\varepsilon}}:(0,\varepsilon_0)\to(0,\infty),$
${\varepsilon}\in(0,{\varepsilon}_0),$
$\varepsilon_0=\varepsilon(z_0)>0,$ such that
\begin{equation}\label{eqI}
I_{z_0}({\varepsilon})\ \colon =\
\int\limits_{{\varepsilon}}^{{\varepsilon}_0}{\psi}_{z_0,{\varepsilon}}(t)\
dt\ <\ \infty\ \ \ \ \ \ \forall\
{\varepsilon}\in(0,{\varepsilon}_0)\end{equation} and
\begin{equation}\label{eqT}
\int\limits_{{\varepsilon}<|z-z_0|<{\varepsilon}_0}\
K^T_{{\mu}}(z,z_0)\cdot{\psi}^2_{z_0,{\varepsilon}}(|z-z_0|)\ dm(z)\
=\ o(I^2_{z_0}({\varepsilon})) \ \ \  \ \ \ \hbox{ as
${\varepsilon}\to 0$}\ \ \ \forall\ z_0\in S\ .\end{equation} Then
the Beltrami equation (\ref{eqBeltrami}) has a regular homeomorphic
solution $f_{\mu}$ with the hydrodynamic normalization
$f_{\mu}(z)=z+o(1)$ as $z\to\infty$.}

\bigskip

Here and further $dm(z)=dxdy$, $z=x+iy$, corresponds to the Lebesgue
measure in ${\Bbb C}$.

\bigskip

\begin{proof}
By Lemma 3 and Remark 2 in \cite{RSY$_4$} the Beltrami equation
(\ref{eqBeltrami}) has a regular homeomorphic solution $f$ in
$\mathbb C$ under the hypotheses on $\mu$ given above. Note that $f$
is holomorphic and univalent (one-to-one), i.e. conformal, and with
no zeros outside of a closed disk $$\{ z\in\mathbb C: |z|\leq R\}, \
\ \ \ R>0\ ,$$ because the support $S$ of $\mu$ is compact.

Let us consider the function $F(\zeta):=f(1/\zeta)$,
$\zeta\in\mathbb C_0:=\overline{\mathbb C}\setminus\{ 0\}$,
$\overline{\mathbb C}=\mathbb C\cup\{\infty\}$, that is conformal in
a punctured disk $\mathbb D_r\setminus\{ 0\}$, where $\mathbb
D_r=\{\zeta\in\mathbb C: |\zeta|<r\}$, $r=1/R$, and $0$ is its
isolated singular point. In view of the Casorati-Weierstrass
theorem, see e.g. Proposition II.6.3 in \cite{FL}, $0$ cannot be
essential singular point because the mapping $F$ is homeomorphic.

Moreover, $0$ cannot be a removable singular point of $F$. Indeed,
let us assume that $F$ has a finite limit $\lim\limits_{\zeta\to
0}F(\zeta)=c$. Then the extended mapping $\tilde F$ is a
homeomorphism of $\overline{\mathbb C}$ into $\mathbb C$. However,
by stereographic projection $\overline{\mathbb C}$ is homeomorphic
to the sphere $\mathbb S^2$ and, consequently, by the Brouwer
theorem on the invariance of domain the set $C:=\tilde
F(\overline{\mathbb C})$ is open in $\overline{\mathbb C}$, see e.g.
Theorem 4.8.16 in \cite{Sp}. In addition, the set $C$ is compact as
a continuous image of the compact space $\overline{\mathbb C}$.
Hence the set $\overline{\mathbb C}\setminus C\ne\emptyset$ is also
open in $\overline{\mathbb C}$. The latter contradicts the
connectivity of $\overline{\mathbb C}$, see e.g. Proposition I.1.1
in \cite{FL}.

Thus, $0$ is a (unique) pole of the function $F$ in the disk
$\mathbb D_r$. Hence the function $\Phi(\zeta):=1/F(\zeta)$ has a
removable singularity at $0$ and $\Phi(0)=0$. By the Riemann
extension theorem, see e.g. Proposition II.3.7 in \cite{FL}, the
extended function $\tilde{\Phi}$ is conformal in $\mathbb D_r$. By
the Rouche theorem $\tilde{\Phi}^{\prime}(0)\ne 0$, see e.g. Theorem
63 in \cite{TV}, and, consequently, the function $\tilde{\Phi}$ has
the expansion of the form $c_1\zeta+c_2\zeta^2+\ldots$ in the disk
$\mathbb D_r$ with $c_1\ne 0$. Consequently, along the set $\{
z\in\mathbb C: |z|>R\}$
$$
f(z)\ =\ \frac{1}{\Phi(\frac{1}{z})}\ =
\frac{1}{c_1z^{-1}+c_2z^{-2}+\ldots}\ =\
\frac{z}{c_1}\left(1+\frac{c_2}{c_1}z^{-1}+\ldots\right)^{-1}\ =\
c_1^{-1}z\ -\ c_1^{-2}c_2\ +\ o(1)\ ,
$$
i.e. the function $f_{\mu}(z):=c_1f(z)+c_2/c_1$ gives the desired
regular homeomorphic solution of the Beltrami equation with the
hydrodynamic normalization $f_{\mu}(z)=z+o(1)$ as $z\to\infty$.
\end{proof} $\Box$

\bigskip

In particular, by relations (\ref{eqConnect}) we obtain from Lemma 1
the following consequence.

\bigskip

{\bf Corollary 1.} {\it Let a function $\mu : \mathbb C\to{\Bbb C}$
be with compact support $S$, $|\mu (z)| < 1$  a.e. and $K_{\mu}\in
L^1(S)$ and let ${\psi}:(0,\varepsilon_0)\to(0,\infty)$ for some
$\varepsilon_0>0$ be a measurable function such that
\begin{equation}\label{eqInt} \int\limits_{0}^{\varepsilon_0}{\psi}(t)\ dt
 =\ \infty\ ,\ \ \ \ \ \
\ \ \int\limits_{\varepsilon}^{\varepsilon_0}{\psi}(t)\ dt\ <\
\infty\ \ \ \ \ \ \ \ \ \forall\ \varepsilon\in(0,\varepsilon_0)\ .
\end{equation} Suppose that \begin{equation}\label{eqK} \int\limits_{{\varepsilon}<|z-z_0|<{\varepsilon}_0}
K_{{\mu}}(z)\cdot{\psi}^2(|z-z_0|)\ dm(z)\ \le\ O\left(
\int\limits_{{\varepsilon}}^{{\varepsilon}_0}\ {\psi}(t)\ dt\right)\
\ \ \ \ \hbox{as ${\varepsilon}\to 0$} \ \ \ \forall\ z_0\in S\
.\end{equation} Then the Beltrami equation (\ref{eqBeltrami}) has a
regular homeomorphic solution $f$ with the hydrodynamic
normalization $f(z)=z+o(1)$ as $z\to\infty$.}

\bigskip

Recall that a real-valued function $u$ in a domain $D$ in ${\Bbb C}$
is said to be of {\bf bounded mean oscillation} in $D$, abbr.
$u\in{\rm BMO}(D)$, if $u\in L_{\rm loc}^1(D)$ and
\begin{equation}\label{lasibm_2.2_1}\Vert u\Vert_{*}:=
\sup\limits_{B}{\frac{1}{|B|}}\int\limits_{B}|u(z)-u_{B}|\,dm(z)<\infty\,,\end{equation}
where the supremum is taken over all discs $B$ in $D$ and
$$u_{B}={\frac{1}{|B|}}\int\limits_{B}u(z)\,dm(z)\,.$$ We write $u\in{\rm BMO}_{\rm loc}(D)$ if
$u\in{\rm BMO}(U)$ for every relatively compact subdomain $U$ of $D$
(we also write BMO or ${\rm BMO}_{\rm loc }$ if it is clear from the
context what $D$ is).

The class BMO was introduced by John and Nirenberg (1961) in the
paper \cite{JN} and soon became an important concept in harmonic
analysis, partial differential equations and related areas, see e.g.
\cite{HKM} and \cite{RR}.
\medskip

A function $\varphi$ in BMO is said to have {\bf vanishing mean
oscillation}, abbr. $\varphi\in{\rm VMO}$, if the supremum in
(\ref{lasibm_2.2_1}) taken over all balls $B$ in $D$ with
$|B|<\varepsilon$ converges to $0$ as $\varepsilon\to0$. VMO has
been introduced by Sarason in \cite{Sar}. There are a number of
papers devoted to the study of partial differential equations with
coefficients of the class VMO, see e.g. \cite{CFL}, \cite{IS},
\cite{MRV}, \cite{Pal}, \cite{Ra$_1$} and \cite{Ra$_2$}.

\medskip

{\bf Remark 1.} Note that $W^{\,1,2}\left({{D}}\right) \subset VMO
\left({{D}}\right),$ see e.g. \cite{BN}.

\medskip

Following \cite{IR}, we say that a function $\varphi:D\to{\Bbb R}$
has {\bf finite mean oscillation} at a point $z_0\in D$, abbr.
$\varphi\in{\rm FMO}(z_0)$, if
\begin{equation}\label{FMO_eq2.4}\overline{\lim\limits_{\varepsilon\to0}}\ \ \
\dashint_{B(z_0,\varepsilon)}|{\varphi}(z)-\widetilde{\varphi}_{\varepsilon}(z_0)|\,dm(z)<\infty\,,\end{equation}
where \begin{equation}\label{FMO_eq2.5}
\widetilde{\varphi}_{\varepsilon}(z_0)=\dashint_{B(z_0,\varepsilon)}
{\varphi}(z)\,dm(z)\end{equation} is the mean value of the function
${\varphi}(z)$ over the disk $B(z_0,\varepsilon):=\{ z\in\mathbb C:
|z-z_0|<\varepsilon\}$. Note that the condition (\ref{FMO_eq2.4})
includes the assumption that $\varphi$ is integrable in some
neighborhood of the point $z_0$. We say also that a function
$\varphi:D\to{\Bbb R}$ is of {\bf finite mean oscillation in $D$},
abbr. $\varphi\in{\rm FMO}(D)$ or simply $\varphi\in{\rm FMO}$, if
$\varphi\in{\rm FMO}(z_0)$ for all points $z_0\in D$. We write
$\varphi\in{\rm FMO}(\overline{D})$ if $\varphi$ is given in a
domain $G$ in $\Bbb{C}$ such that $\overline{D}\subset G$ and
$\varphi\in{\rm FMO}(G)$.

\medskip

The following statement is obvious by the triangle inequality.

\medskip

{\bf Proposition 1.} {\it If, for a  collection of numbers
$\varphi_{\varepsilon}\in{\Bbb R}$,
$\varepsilon\in(0,\varepsilon_0]$,
\begin{equation}\label{FMO_eq2.7}\overline{\lim\limits_{\varepsilon\to0}}\ \ \
\dashint_{B(z_0,\varepsilon)}|\varphi(z)-\varphi_{\varepsilon}|\,dm(z)<\infty\,,\end{equation}
then $\varphi $ is of finite mean oscillation at $z_0$.}

\medskip

In particular, choosing here  $\varphi_{\varepsilon}\equiv0$,
$\varepsilon\in(0,\varepsilon_0]$ in Proposition 1, we obtain the
following.

\medskip

{\bf Corollary 2.} {\it If, for a point $z_0\in D$,
\begin{equation}\label{FMO_eq2.8}\overline{\lim\limits_{\varepsilon\to 0}}\ \ \
\dashint_{B(z_0,\varepsilon)}|\varphi(z)|\,dm(z)<\infty\,,
\end{equation} then $\varphi$ has finite mean oscillation at
$z_0$.}

\medskip

Recall that a point $z_0\in D$ is called a {\bf Lebesgue point} of a
function $\varphi:D\to{\Bbb R}$ if $\varphi$ is integrable in a
neighborhood of $z_0$ and \begin{equation}\label{FMO_eq2.7a}
\lim\limits_{\varepsilon\to 0}\ \ \ \dashint_{B(z_0,\varepsilon)}
|\varphi(z)-\varphi(z_0)|\,dm(z)=0\,.\end{equation} It is known
that, almost every point in $D$ is a Lebesgue point for every
function $\varphi\in L^1(D)$. Thus, we have by Proposition 1 the
next corollary.

\medskip

{\bf Corollary 3.} {\it Every locally integrable function
$\varphi:D\to{\Bbb R}$ has a finite mean oscillation at almost every
point in $D$.}

\medskip

{\bf Remark 2.} Note that the function
$\varphi(z)=\log\left(1/|z|\right)$ belongs to BMO in the unit disk
$\Delta$, see, e.g., \cite{RR}, p. 5, and hence also to FMO.
However, $\widetilde{\varphi}_{\varepsilon}(0)\to\infty$ as
$\varepsilon\to0$, showing that condition (\ref{FMO_eq2.8}) is only
sufficient but not necessary for a function $\varphi$ to be of
finite mean oscillation at $z_0$. Clearly, ${\rm BMO}(D)\subset{\rm
BMO}_{\rm loc}(D)\subset{\rm FMO}(D)$ and as well-known ${\rm
BMO}_{\rm loc}\subset L_{\rm loc}^p$ for all $p\in[1,\infty)$, see,
e.g., \cite{JN} or \cite{RR}. However, FMO is not a subclass of
$L_{\rm loc}^p$ for any $p>1$ but only of $L_{\rm loc}^1$. Thus, the
class FMO is much more wider than ${\rm BMO}_{\rm loc}$.

\medskip

Versions of the next lemma has been first proved for the class BMO
in \cite{RSY}. For the FMO case, see the papers \cite{IR, RS,
RSY$_2$, RSY$_3$}  and the monographs \cite{GRSY} and \cite{MRSY}.

\medskip

{\bf Lemma 2.} {\it Let $D$ be a domain in ${\Bbb C}$ and let
$\varphi:D\to{\Bbb R}$ be a  non-negative function  of the class
${\rm FMO}(z_0)$ for some $z_0\in D$. Then
\begin{equation}\label{eq13.4.5}\int\limits_{\varepsilon<|z-z_0|<\varepsilon_0}\frac{\varphi(z)\,dm(z)}
{\left(|z-z_0|\log\frac{1}{|z-z_0|}\right)^2}=O\left(\log\log\frac{1}{\varepsilon}\right)\
\quad\text{as}\quad\varepsilon\to 0\end{equation} for some
$\varepsilon_0\in(0,\delta_0)$ where $\delta_0=\min(e^{-e},d_0)$,
$d_0=\sup\limits_{z\in D}|z-z_0|$.}

\medskip

The following statement will be also useful later on, see e.g.
Theorem 3.2 in \cite{RSY$_5$}.

\medskip

{\bf Proposition 2.} {\it Let $Q:{\Bbb D}\to[0,\infty]$ be a
measurable function such that
\begin{equation}\label{eq5.555} \int\limits_{\Bbb
D}\Phi(Q(z))\,dm(z)<\infty\end{equation} where
$\Phi:[0,\infty]\to[0,\infty]$ is a non-decreasing convex function
such that \begin{equation}\label{eq3.333a}
\int\limits_{\delta}^{\infty}\frac{d\tau}{\tau\Phi^{-1}(\tau)}=\infty\end{equation}
for some $\delta>\Phi(+0)$. Then \begin{equation}\label{eq3.333A}
\int\limits_0^1\frac{dr}{rq(r)}=\infty\end{equation} where $q(r)$ is
the average of the function $Q(z)$ over the circle $|z|=r$.}

\medskip

Here we use the following notions of the inverse function for
monotone functions. Namely, for every non-decreasing function
$\Phi:[0,\infty]\to[0,\infty]$ the inverse function
$\Phi^{-1}:[0,\infty]\to[0,\infty]$ can be well-defined by setting
\begin{equation}\label{eqINVERSE}
\Phi^{-1}(\tau)\ :=\ \inf\limits_{\Phi(t)\geq\tau} t
\end{equation}
Here $\inf$ is equal to $\infty$ if the set of $t\in[0,\infty]$ such
that $\Phi(t)\geq\tau$ is empty. Note that the function $\Phi^{-1}$
is non-decreasing, too. It is evident immediately by the definition
that $\Phi^{-1}(\Phi(t)) \leq t$ for all $t\in[0,\infty]$ with the
equality except intervals of constancy of the function $\Phi(t)$.

\bigskip

Finally, recall connections between some integral conditions, see
e.g. Theorem 2.5 in \cite{RSY$_5$}.

\medskip

{\bf Remark 3.} Let $\Phi:[0,\infty]\to[0,\infty]$ be a
non-decreasing function and set
\begin{equation}\label{eqLOGFi}
H(t)\ =\ \log\Phi(t)\ .
\end{equation}
Then the equality
\begin{equation}\label{eq333Frer}\int\limits_{\Delta}^{\infty}H'(t)\,\frac{dt}{t}=\infty,
\end{equation}
implies the equality
\begin{equation}\label{eq333F}\int\limits_{\Delta}^{\infty}
\frac{dH(t)}{t}=\infty\,,\end{equation} and (\ref{eq333F}) is
equivalent to
\begin{equation}\label{eq333B}
\int\limits_{\Delta}^{\infty}H(t)\,\frac{dt}{t^2}=\infty\,\end{equation}
for some $\Delta>0$, and (\ref{eq333B}) is equivalent to each of the
equalities
\begin{equation}\label{eq333C} \int\limits_{0}^{\delta_*}H\left(\frac{1}{t}\right)\,{dt}=\infty\end{equation} for
some $\delta_*>0$, \begin{equation}\label{eq333D}
\int\limits_{\Delta_*}^{\infty}\frac{d\eta}{H^{-1}(\eta)}=\infty\end{equation}
for some $\Delta_*>H(+0)$ and to (\ref{eq3.333a}) for some
$\delta>\Phi(+0)$.

Moreover, (\ref{eq333Frer}) is equivalent to (\ref{eq333F}) and to
hence (\ref{eq333Frer})–(\ref{eq333D}) as well as to
(\ref{eq3.333a}) are equivalent to each other if $\Phi$ is in
addition absolutely continuous. In particular, all the given
conditions are equivalent if $\Phi$ is convex and non-decreasing.

%\medskip

Note that the integral in (\ref{eq333F}) is understood as the
Lebesgue--Stieltjes integral and the integrals in (\ref{eq333Frer})
and (\ref{eq333B})--(\ref{eq333D}) as the ordinary Lebesgue
integrals. It is necessary to give one more explanation. From the
right hand sides in the conditions (\ref{eq333Frer})--(\ref{eq333D})
we have in mind $+\infty$. If $\Phi(t)=0$ for $t\in[0,t_*$, then
$H(t)=-\infty$ for $t\in[0,t_*]$ and we complete the definition
$H'(t)=0$ for $t\in[0,t_*]$. Note, the conditions (\ref{eq333F}) and
(\ref{eq333B}) exclude that $t_*$ belongs to the interval of
integrability because in the contrary case the left hand sides in
(\ref{eq333F}) and (\ref{eq333B}) are either equal to $-\infty$ or
indeterminate. Hence we may assume in
(\ref{eq333Frer})--(\ref{eq333C}) that $\delta >t_0$,
correspondingly, $\Delta<1/t_0$ where
$t_0:=\sup\limits_{\Phi(t)=0}t$, and set $t_0=0$ if $\Phi(0)>0$.

%\medskip

The most interesting of the above conditions is (\ref{eq333B}) that
can be rewritten in the form:
\begin{equation}\label{eq5!}
\int\limits_{\Delta}^{\infty}\log\, \Phi(t)\ \frac{dt}{t^{2}}\ =\
+\infty\ \ \ \ \ \ \mbox{for some $\Delta > 0$}\ .
\end{equation}

\section{The Dirichlet problem in simply connected domains}

Recall that a mapping $f:D\to{\Bbb C}$ is called {\bf discrete} if
the preimage  $f^{-1}(y)$ consists of isolated points for every
$y\in{\Bbb C}$, and {\bf open} if $f$ maps every open set
$U\subseteq D$ onto an open set in ${\Bbb C}$. If
$\varphi(\zeta)\not\equiv{\rm const}$, then the {\bf regular
solution}  of the Dirichlet problem (\ref{eqDIRICHLET}) for the
Beltrami equation (\ref{eqBeltrami}) is a continuous, discrete and
open mapping $f:D\to{\Bbb C}$ of the Sobolev class $W_{\rm
loc}^{1,1}$ with its Jacobian $J_f(z)=|f_z|^2-|f_{\bar z}|^2\neq0$
a.e. satisfying (\ref{eqBeltrami}) a.e. and the condition
(\ref{eqDIRICHLET}). The regular solution of such a problem with
$\varphi(\zeta)\equiv c$, $\zeta\in\partial D$, for the Beltrami
equation (\ref{eqBeltrami}) is the function $f(z)\equiv c$, $z\in
D$.

In this section, we prove that a regular solution of the Dirichlet
problem (\ref{eqDIRICHLET}) exists for every continuous function
$\varphi:\partial D\to{\Bbb R}$ for wide classes of the degenerate
Beltrami equations (\ref{eqBeltrami}) in an arbitrary bounded domain
$D$ with no boundary component degenerated to a point, and that such
a solution can be represented in the form of the composition of a
regular homeomorphic solution of (\ref{eqBeltrami}) with
hydrodynamic normalization and a holomorphic solution of the
Dirichlet problem associated with it. The main criteria are
formulated by us in terms of the tangent dilatations
$K^T_{\mu}(z,z_0)$ which are more refined although the corresponding
criteria remain valid for the usual dilatation $K_{\mu}(z)$.

\medskip

We assume further that the dilatations $K^T_{\mu}(z,z_0)$ and
$K_{\mu}(z)$ are extended by $1$ outside of the domain $D$.

\medskip

{\bf Lemma 4.} {\it Let $D$ be a bounded simply connected domain in
${\Bbb C}.$ Suppose that $\mu:D\to{\Bbb C}$ is a measurable function
with $|\mu(z)|<1$ a.e., $K_{\mu}\in L^1(D)$ and
\begin{equation}\label{3omal}
\int\limits_{\varepsilon<|z-z_0|<\varepsilon_0}
K^T_{\mu}(z,z_0)\cdot\psi^2_{z_0,\varepsilon}(|z-z_0|)\,dm(z)=o(I_{z_0}^{2}(\varepsilon))\quad{\rm
as}\quad\varepsilon\to0\ \ \forall\ z_0\in\overline{D}\end{equation}
for some $\varepsilon_0=\varepsilon(z_0)>0$ and a family of
measurable functions $\psi_{z_0,\varepsilon}:
(0,\varepsilon_0)\to(0,\infty)$ with
\begin{equation}\label{eq3.5.3}
I_{z_0}(\varepsilon)\colon
=\int\limits_{\varepsilon}^{\varepsilon_0}
\psi_{z_0,\varepsilon}(t)\,dt<\infty\qquad\forall\
\varepsilon\in(0,\varepsilon_0)\,.\end{equation} Then the Beltrami
equation (\ref{eqBeltrami}) has a regular solution $f$ of the
Dirichlet problem (\ref{eqDIRICHLET}) for each continuous function
$\varphi:\partial D\to{\Bbb R}$.

Moreover, such a solution $f$ can be represented as the composition
\begin{equation}\label{eqHYDRO}
f\ =\ h\circ g\ ,\ \ \ \ \ \ \hbox{$g(z)\, =\ z\, +\, o(1)$\ \ \ as\
\ \ $z\to\infty$}\ ,
\end{equation}
where $g:\mathbb C\to\mathbb C$ is a regular homeomorphic solution
of the Beltrami equation (\ref{eqBeltrami}) in $\mathbb C$ with
$\mu$ extended by zero outside of $D$ and $h:D_*\to\mathbb C,$
$D_*:=g(D)$, is a holomorphic solution of the Dirichlet problem
\begin{equation}\label{eqHOLOMORPHIC}
\lim_{\xi\to\zeta}\ {\rm Re}\, h(\xi)\ =\ \varphi_*(\zeta)\ \ \ \ \
\forall\ \zeta\in\partial D_*\ ,\ \ \ \ \ \mbox{where
$\varphi_*:=\varphi\circ g^{-1}$.}
\end{equation}}

\bigskip

\begin{proof} Indeed, by Lemma 1 there is a regular homeomorphic
solution with hydrodynamic normalization $g(z):=z + o(1)$ as
$z\to\infty$ of the Beltrami equation (\ref{eqBeltrami}) in $\mathbb
C$ with $\mu$ extended by zero outside of $D$. Note that $D_*:=g(D)$
is also a simply connected domain in $\mathbb C$ with no boundary
component degenerated to a single point because of $g:\mathbb
C\to\mathbb C$ is a homeomorphism. Consequently, by Theorem 4.2.2
and Corollary 4.1.8 in \cite{Rans} there is a unique harmonic
function $u:D_*\to\mathbb R$ that satisfies the Dirichlet boundary
condition
\begin{equation}\label{eqHARMONIC}
\lim_{\xi\to\zeta}\ u(\xi)\ :=\ \varphi_*(\zeta)\ \ \ \ \ \forall\
\zeta\in\partial D_*\ ,\ \ \ \ \ \mbox{where
$\varphi_*:=\varphi\circ g^{-1}$.}
\end{equation}
On the other hand, there is a conjugate harmonic function
$v:D_*\to\mathbb R$ such that $h:=u+iv:D_*\to\mathbb C$ forms a
holomorphic function because of the domain $D_*$ is simply
connected, see e.g. arguments in the beginning of the book
\cite{Ko}. Thus, the function $f:=h\circ g$ gives the desired
solution of the Dirichlet problem (\ref{eqDIRICHLET}) for the
Beltrami equation (\ref{eqBeltrami}).
\end{proof} $\Box$

\medskip

{\bf Remark 4.} Note that if the family of the functions
$\psi_{z_0,\varepsilon}(t)\equiv\psi_{z_0}(t)$ is independent on the
parameter $\varepsilon$, then the condition (\ref{3omal}) implies
that $I_{z_0}(\varepsilon)\to \infty$ as $\varepsilon\to 0$. This
follows immediately from arguments by contradiction, apply for it
(\ref{eqConnect}) and the condition $K_{\mu}\in L^1(D)$. Note also
that (\ref{3omal}) holds, in particular, if, for some
$\varepsilon_0=\varepsilon(z_0)$,
\begin{equation}\label{333omal}
\int\limits_{|z-z_0|<\varepsilon_0}
K^T_{\mu}(z,z_0)\cdot\psi_{z_0}^2(|z-z_0|)\,dm(z)<\infty \qquad
\forall\ z_0\in\overline{D}\end{equation} and
$I_{z_0}(\varepsilon)\to \infty$ as $\varepsilon\to 0$. In other
words, for the solvability of the Dirichlet problem
(\ref{eqDIRICHLET}) for the Beltrami equation (\ref{eqBeltrami}) for
all continuous boundary functions $\varphi$, it is sufficient that
the integral in (\ref{333omal}) converges for some nonnegative
function $\psi_{z_0}(t)$ that is locally integrable over
$(0,\varepsilon_0 ]$ but has a nonintegrable singularity at $0$. The
functions $\log^{\lambda}(e/|z-z_0|)$, $\lambda\in (0,1)$, $z\in\Bbb
D$, $z_0\in\overline{\Bbb D}$, and $\psi(t)=1/(t \,\, \log(e/t))$,
$t\in(0,1)$, show that the condition (\ref{333omal}) is compatible
with the condition $I_{z_0}(\varepsilon)\to\infty $ as
$\varepsilon\to 0$. Furthermore, the condition (\ref{3omal}) shows
that it is sufficient for the solvability of the Dirichlet problem
even if the integral in (\ref{333omal}) is divergent in a controlled
way.

\medskip

Choosing $\psi(t)=1/\left(t\, \log\left(1/t\right)\right)$ in Lemma
4, we obtain by Lemma 2 the following result.

\medskip

{\bf Theorem 1.} {\it Let $D$ be a bounded simply connected domain
in $\mathbb C$ and $\mu:D\to{\Bbb C}$ be a measurable function with
$|\mu(z)|<1$ a.e. and $K_{\mu}\in L^1(D)$. Suppose that
$K^T_{\mu}(z,z_0)\leqslant Q_{z_0}(z)$ a.e. in $U_{z_0}$ for every
point $z_0\in \overline{D}$, a neighborhood $U_{z_0}$ of $z_0$ and a
function $Q_{z_0}: U_{z_0}\to[0,\infty]$ in the class ${\rm
FMO}({z_0})$. Then the Beltrami equation (\ref{eqBeltrami}) has a
regular solution of the Dirichlet problem (\ref{eqDIRICHLET}) with
the representation (\ref{eqHYDRO}) for each continuous function
$\varphi:\partial D\to{\Bbb R}$.}

\medskip

In particular, by Proposition 1 the conclusion of Theorem 1 holds if
every point $z_0\in\overline{D}$ is the Lebesgue point of the
function $Q_{z_0}$.

\bigskip

By Corollary 2 we obtain the following nice consequence of Theorem
1, too.

\medskip

{\bf Corollary 4.} {\it Let $D$ be a bounded simply connected domain
in $\mathbb C$ and $\mu:D\to{\Bbb C}$ be a measurable function with
$|\mu(z)|<1$ a.e., $K_{\mu}\in L^1(D)$ and
\begin{equation}\label{eqMEAN}\overline{\lim\limits_{\varepsilon\to0}}\quad
\dashint_{B(z_0,\varepsilon)}K^T_{\mu}(z,z_0)\,dm(z)<\infty\qquad\forall\
z_0\in\overline{D}\, .\end{equation}Then the Beltrami equation
(\ref{eqBeltrami}) has a regular solution of the Dirichlet problem
(\ref{eqDIRICHLET}) with the representation (\ref{eqHYDRO}) for each
continuous function $\varphi:\partial D\to{\Bbb R}$.}

\medskip

Since $K^T_{\mu}(z,z_0) \leqslant K_{\mu}(z)$ for all $z$ and
$z_0\in \Bbb C$, we also obtain the following consequences of
Theorem 1.

\medskip

{\bf Corollary 5.} {\it Let $D$ be a bounded simply connected domain
in $\mathbb C$ and $\mu:D\to{\Bbb C}$ be a measurable function with
$|\mu(z)|<1$ a.e. and $K_{\mu}$ have a dominant $Q:\mathbb
C\to[1,\infty)$ in the class {\rm BMO}$_{\rm loc}$. Then the
Beltrami equation (\ref{eqBeltrami}) has a regular solution of the
Dirichlet problem (\ref{eqDIRICHLET}) with the representation
(\ref{eqHYDRO}) for each continuous function $\varphi:\partial
D\to{\Bbb R}$.}

\medskip

{\bf Remark 5.} In particular, the conclusion of Corollary 5 holds
if $Q\in{\rm W}^{1,2}_{\rm loc}$ because $W^{\,1,2}_{\rm loc}
\subset {\rm VMO}_{\rm loc}$, see e.g. \cite{BN}.

\medskip

{\bf Corollary 6.} {\it Let $D$ be a bounded simply connected domain
in $\mathbb C$ and $\mu:D\to{\Bbb C}$ be a measurable function with
$|\mu(z)|<1$ and $K_{\mu}(z)\leqslant Q(z)$ a.e. in $D$ with a
function $Q$ in the class ${\rm FMO}(\overline{D})$. Then the
Beltrami equation (\ref{eqBeltrami}) has a regular solution of the
Dirichlet problem (\ref{eqDIRICHLET}) with the representation
(\ref{eqHYDRO}) for each continuous function $\varphi:\partial
D\to{\Bbb R}$.}

\medskip

Similarly, choosing in Lemma 4 the function $\psi(t)=1/t$, we come
to the next statement.

\medskip

{\bf Theorem 2.} {\it Let $D$ be a bounded simply connected domain
in $\mathbb C$ and $\mu:D\to{\Bbb C}$ be a measurable function with
$|\mu(z)|<1$ a.e. and $K_{\mu}\in L^1(D)$. Suppose that
\begin{equation}\label{eqLOG}
\int\limits_{\varepsilon<|z-z_0|<\varepsilon_0}K^T_{\mu}(z,z_0)\,\frac{dm(z)}{|z-z_0|^2}
=o\left(\left[\log\frac{1}{\varepsilon}\right]^2\right)\qquad\hbox{as
$\varepsilon\to 0$}\qquad\forall\ z_0\in\overline{D}\end{equation}
for some $\varepsilon_0=\varepsilon(z_0)>0$. Then Beltrami equation
(\ref{eqBeltrami}) has a regular solution of the Dirichlet problem
(\ref{eqDIRICHLET}) with the representation (\ref{eqHYDRO}) for each
continuous function $\varphi:\partial D\to{\Bbb R}$.}

\medskip

{\bf Remark 6.} Choosing in Lemma 4 the function
$\psi(t)=1/(t\log{1/t})$ instead of $\psi(t)=1/t$, we are able to
replace (\ref{eqLOG}) by
\begin{equation}\label{eqLOGLOG}
\int\limits_{\varepsilon<|z-z_0|<\varepsilon_0}\frac{K^T_{\mu}(z,z_0)\,dm(z)}
{\left(|z-z_0|\log{\frac{1}{|z-z_0|}}\right)^2}
=o\left(\left[\log\log\frac{1}{\varepsilon}\right]^2\right)\end{equation}
In general, we are able to give here the whole scale of the
corresponding conditions in $\log$ using functions $\psi(t)$ of the
form
$1/(t\log{1}/{t}\cdot\log\log{1}/{t}\cdot\ldots\cdot\log\ldots\log{1}/{t})$.

\medskip

Choosing in Lemma 4 the functional parameter
${\psi}_{z_0,{\varepsilon}}(t) \equiv {\psi}_{z_0}(t) \colon =
1/[tk^T_{\mu}(z_0,t)]$, where $k_{\mu}^T(z_0,r)$ is the integral
mean of $K^T_{{\mu}}(z,z_0)$ over the circle $S(z_0,r)\, :=\, \{ z
\in\mathbb C:\, |z-z_0|\,=\, r\}$, we obtain one more important
conclusion.

\medskip

{\bf Theorem 3.} {\it Let $D$ be a bounded simply connected domain
in $\mathbb C$ and $\mu:D\to{\Bbb C}$ be a measurable function with
$|\mu(z)|<1$ a.e. and $K_{\mu}\in L^1(D)$. Suppose that
\begin{equation}\label{eqLEHTO}\int\limits_{0}^{\varepsilon_0}
\frac{dr}{rk^T_{\mu}(z_0,r)}=\infty\qquad\forall\
z_0\in\overline{D}\end{equation} for some
$\varepsilon_0=\varepsilon(z_0)>0$. Then Beltrami equation
(\ref{eqBeltrami}) has a regular solution of the Dirichlet problem
(\ref{eqDIRICHLET}) with the representation (\ref{eqHYDRO}) for each
continuous function $\varphi:\partial D\to{\Bbb R}$.}

\medskip

{\bf Corollary 7.} {\it Let $D$ be a bounded simply connected domain
in $\mathbb C$ and $\mu:D\to{\Bbb C}$ be a measurable function with
$|\mu(z)|<1$ a.e., $K_{\mu}\in L^1(D)$ and
\begin{equation}\label{eqLOGk}k^T_{\mu}(z_0,\varepsilon)=O\left(\log\frac{1}{\varepsilon}\right)
\qquad\mbox{as}\ \varepsilon\to0\qquad\forall\ z_0\in\overline{D}\
.\end{equation} Then the Beltrami equation (\ref{eqBeltrami}) has a
regular solution of the Dirichlet problem (\ref{eqDIRICHLET}) with
the representation (\ref{eqHYDRO}) for each continuous function
$\varphi:\partial D\to{\Bbb R}$.}

\medskip

{\bf Remark 7.} In particular, the conclusion of Corollary 7 holds
if
\begin{equation}\label{eqLOGK} K^T_{\mu}(z,z_0)=O\left(\log\frac{1}{|z-z_0|}\right)\qquad{\rm
as}\quad z\to z_0\quad\forall\ z_0\in\overline{D}\,.\end{equation}
Moreover, the condition (\ref{eqLOGk}) can be replaced by the whole
series of more weak conditions
\begin{equation}\label{edLOGLOGk}
k^T_{\mu}(z_0,\varepsilon)=O\left(\left[\log\frac{1}{\varepsilon}\cdot\log\log\frac{1}
{\varepsilon}\cdot\ldots\cdot\log\ldots\log\frac{1}{\varepsilon}
\right]\right) \qquad\forall\ z_0\in \overline{D}\ .
\end{equation}

\medskip

Combining Theorems 3, Proposition 2 and Remark 3, we obtain the
following result.

\medskip

{\bf Theorem 4.} {\it Let $D$ be a bounded simply connected domain
in $\mathbb C$ and $\mu:D\to{\Bbb C}$ be a measurable function with
$|\mu(z)|<1$ a.e. and $K_{\mu}\in L^1(D)$. Suppose that
\begin{equation}\label{eqINTEGRAL}\int\limits_{U_{z_0}}\Phi_{z_0}\left(K^T_{\mu}(z,z_0)\right)\,dm(z)<\infty
\qquad\forall\ z_0\in \overline{D}\end{equation} for a neighborhood
$U_{z_0}$ of $z_0$ and a convex non-decreasing function
$\Phi_{z_0}:[0,\infty]\to[0,\infty]$ with
\begin{equation}\label{eqINT}
\int\limits_{\Delta(z_0)}^{\infty}\log\,\Phi_{z_0}(t)\,\frac{dt}{t^2}\
=\ +\infty\end{equation} for some $\Delta(z_0)>0$. Then Beltrami
equation (\ref{eqBeltrami}) has a regular solution of the Dirichlet
problem (\ref{eqDIRICHLET}) with the representation (\ref{eqHYDRO})
for each continuous function $\varphi:\partial D\to{\Bbb R}$.}

\medskip

{\bf Corollary 8.} {\it Let $D$ be a bounded simply connected domain
in $\mathbb C$ and $\mu:D\to{\Bbb C}$ be a measurable function with
$|\mu(z)|<1$ a.e., $K_{\mu}\in L^1(D)$ and
\begin{equation}\label{eqEXP}\int\limits_{U_{z_0}}e^{\alpha(z_0) K^T_{\mu}(z,z_0)}\,dm(z)<\infty
\qquad\forall\ z_0\in \overline{D}\end{equation} for some
$\alpha(z_0)>0$ and a neighborhood $U_{z_0}$ of the point $z_0$.
Then the Beltrami equation (\ref{eqBeltrami}) has a regular solution
of the Dirichlet problem (\ref{eqDIRICHLET}) with the representation
(\ref{eqHYDRO}) for each continuous function $\varphi:\partial
D\to{\Bbb R}$.}

\medskip

Since $K^T_{\mu}(z,z_0) \leqslant K_{\mu}(z)$ for $z$ and $z_0\in
\Bbb C$ and $z\in D$, we also obtain the following consequences of
Theorem 4.

\medskip

{\bf Corollary 9.} {\it Let $D$ be a bounded simply connected domain
in $\mathbb C$ and $\mu:D\to{\Bbb C}$ be a measurable function with
$|\mu(z)|<1$ a.e. and $K_{\mu}\in L^1(D)$. Suppose that
\begin{equation}\label{eqINTK}\int\limits_{D}\Phi\left(K_{\mu}(z)\right)\,dm(z)<\infty\end{equation}
for a convex non-decreasing function $\Phi:[0,\infty]\to[0,\infty]$
with
\begin{equation}\label{eqINTF}
\int\limits_{\delta}^{\infty}\log\,\Phi(t)\,\frac{dt}{t^2}\ =\
+\infty\end{equation} for some $\delta>0$. Then Beltrami equation
(\ref{eqBeltrami}) has a regular solution of the Dirichlet problem
(\ref{eqDIRICHLET}) with the representation (\ref{eqHYDRO}) for each
continuous function $\varphi:\partial D\to{\Bbb R}$.}

\medskip

{\bf Corollary 10.} {\it Let $D$ be a bounded simply connected
domain in $\mathbb C$ and $\mu:D\to{\Bbb C}$ be a measurable
function with $|\mu(z)|<1$ a.e. and, for some $\alpha>0$,
\begin{equation}\label{eqEXPA}\int\limits_{D}e^{\alpha K_{\mu}(z)}\,dm(z)\ <\
\infty\ .
\end{equation} Then the Beltrami equation
(\ref{eqBeltrami}) has a regular solution of the Dirichlet problem
(\ref{eqDIRICHLET}) with the representation (\ref{eqHYDRO}) for each
continuous function $\varphi:\partial D\to{\Bbb R}$.}

\medskip

{\bf Remark 8.} By the Stoilow theorem, see e.g. \cite{Sto}, a
regular solution $f$ of the Dirichlet problem (\ref{eqDIRICHLET})
for the Beltrami equation (\ref{eqBeltrami}) with $K_{\mu}\in
L^1_{\rm loc}(D)$ can be represented in the form $f=h\circ F$ where
$h$ is a holomorphic function and $F$ is a homeomorphic regular
solution of (\ref{eqBeltrami}) in the class $W_{\rm loc}^{1,1}$.
Thus, by Theorem 5.1 in \cite{RSY$_5$} the condition (\ref{eqINTF})
is not only sufficient but also necessary to have a regular solution
of the Dirichlet problem (\ref{eqDIRICHLET}) for arbitrary Beltrami
equations (\ref{eqBeltrami}) with the integral constraints
(\ref{eqINTK}) for all continuous functions $\varphi:\partial
D\to\Bbb{R}$, see also Remark 3.

%\medskip

\section{On the Dirichlet problem in general  domains}

In this section we obtain criteria for the existence of multi-valued
solutions $f$ of the Dirichlet problem to the Beltrami equations in
the spirit of the theory of multi-valued analytic functions in
arbitrary bounded domains $D$ in $\mathbb C$ with no boundary
component degenerated to a single point. Our example in Introduction
shows that such domains form the most wide class of domains for
which the problem is always solvable for any continuous boundary
functions.

We say that a discrete open mapping $f:B(z_0,\varepsilon_0)\to{\Bbb
C}$, where $B(z_0,\varepsilon_0)\subseteq D$, is a {\bf local
regular solution of the equation} (\ref{eqBeltrami}) if $f\in W_{\rm
loc}^{1,1}$, $J_f(z)\neq0$ and $f$ satisfies (\ref{eqBeltrami}) a.e.
in $B(z_0,\varepsilon_0)$. The local regular solutions
$f_0:B(z_0,\varepsilon_0)\to{\Bbb C}$ and
$f_*:B(z_*,\varepsilon_*)\to{\Bbb C}$ of the equation
(\ref{eqBeltrami}) will be called extension of each to other if
there is a finite chain of such solutions
$f_i:B(z_i,\varepsilon_i)\to\Bbb{C}$, $i=1,\ldots,m$, such that
$f_1=f_0$, $f_m=f_*$ and $f_i(z)\equiv f_{i+1}(z)$ for $z\in
E_i:=B(z_i,\varepsilon_i)\cap
B(z_{i+1},\varepsilon_{i+1})\neq\emptyset$, $i=1,\ldots,m-1$. A
collection of local regular solutions
$f_j:B(z_j,\varepsilon_j)\to{\Bbb C}$, $j\in J$, will be called a
{\bf multi-valued solution} of the equation (\ref{eqBeltrami}) in
$D$ if the disks $B(z_j,\varepsilon_j)$ cover the whole domain $D$
and $f_j$ are extensions of each to other through the collection and
the collection is maximal by inclusion. A multi-valued solution of
the equation (\ref{eqBeltrami}) will be called a {\bf multi-valued
solution of the Dirichlet problem} (\ref{eqDIRICHLET}) if $u(z)={\rm
Re}\,f(z)={\rm Re}\,f_{j}(z)$, $z\in B(z_j,\varepsilon_j)$, $j\in
J$, is a single-valued function in $D$ satisfying the condition
$\lim\limits_{z\in\zeta}u(z)=\varphi(\zeta)$ for all $\zeta\to
\partial D$.

\medskip

As it was before, we assume further that the dilatations
$K^T_{\mu}(z,z_0)$ and $K_{\mu}(z)$ are extended by $1$ outside of
the domain $D$.

\medskip

{\bf Lemma 5.} {\it Let $D$ be a bounded domain in ${\Bbb C}$ with
no boundary component degenerated to a single point, $\mu:D\to{\Bbb
C}$ be a measurable function with $|\mu(z)|<1$ a.e., $K_{\mu}\in
L^1(D)$ and
\begin{equation}\label{3omalM}
\int\limits_{\varepsilon<|z-z_0|<\varepsilon_0}
K^T_{\mu}(z,z_0)\cdot\psi^2_{z_0,\varepsilon}(|z-z_0|)\,dm(z)=o(I_{z_0}^{2}(\varepsilon))\quad{\rm
as}\quad\varepsilon\to0\ \ \forall\ z_0\in\overline{D}\end{equation}
for some $\varepsilon_0=\varepsilon(z_0)>0$ and a family of
measurable functions $\psi_{z_0,\varepsilon}:
(0,\varepsilon_0)\to(0,\infty)$ with
\begin{equation}\label{eq3.5.3M}
I_{z_0}(\varepsilon)\colon
=\int\limits_{\varepsilon}^{\varepsilon_0}
\psi_{z_0,\varepsilon}(t)\,dt<\infty\qquad\forall\
\varepsilon\in(0,\varepsilon_0)\,.\end{equation} Then the Beltrami
equation (\ref{eqBeltrami}) has a multi-valued solution $f$ of the
Dirichlet problem (\ref{eqDIRICHLET}) for each continuous function
$\varphi:\partial D\to{\Bbb R}$.

Moreover, such a solution $f$ can be represented as the composition
\begin{equation}\label{eqHYDROm}
f\ =\ {\cal A}\circ g\ ,\ \ \ \ \ \ \hbox{$g(z)\, =\ z\, +\, o(1)$\
\ \ as\ \ \ $z\to\infty$}\ ,
\end{equation}
where $g:\mathbb C\to\mathbb C$ is a regular homeomorphic solution
of the Beltrami equation (\ref{eqBeltrami}) in $\mathbb C$ with
$\mu$ extended by zero outside of $D$ and ${\cal A}:D_*\to\mathbb
C,$ $D_*:=g(D)$, is a multi-valued analytic function with a
single-valued harmonic function ${\rm Re}\, {\cal A}$ satisfying the
Dirichlet condition
\begin{equation}\label{eqHOLOMORPHICM}
\lim_{\xi\to\zeta}\ {\rm Re}\, {\cal A}(\xi)\ =\ \varphi_*(\zeta)\ \
\ \ \ \forall\ \zeta\in\partial D_*\ ,\ \ \ \ \ \mbox{where
$\varphi_*:=\varphi\circ g^{-1}$.}
\end{equation}}

\bigskip

\begin{proof} Indeed, by Lemma 1 there is a regular homeomorphic
solution with hydrodynamic normalization $g(z):=z + o(1)$ as
$z\to\infty$ of the Beltrami equation (\ref{eqBeltrami}) in $\mathbb
C$ with $\mu$ extended by zero outside of $D$. Note that $D_*:=g(D)$
is also a simply connected domain in $\mathbb C$ with no boundary
component degenerated to a single point because of $g:\mathbb
C\to\mathbb C$ is a homeomorphism. Consequently, by Theorem 4.2.2
and Corollary 4.1.8 in \cite{Rans} there is a unique harmonic
function $u:D_*\to\mathbb R$ that satisfies the Dirichlet boundary
condition
\begin{equation}\label{eqHARMONICm}
\lim_{\xi\to\zeta}\ u(\xi)\ :=\ \varphi_*(\zeta)\ \ \ \ \ \forall\
\zeta\in\partial D_*\ ,\ \ \ \ \ \mbox{where
$\varphi_*:=\varphi\circ g^{-1}$.}
\end{equation}

Let $B_0=B(z_0,r_0)$ is a disk in the domain $D$. Then ${\frak B}_0
= g(B_0)$ is a simply connected subdomain of the domain $D_*:=g(D)$
where there is a conjugate function $v$ determined up to an additive
constant such that $h=u+iv$ is a single--valued analytic function.
Let us denote through $h_0$ the holomorphic function corresponding
to the choice of such a harmonic function $v_0$ in ${\frak B}_0$
with the normalization $v_0(g(z_0))=0$. Thereby we have determined
the initial element of a multi-valued analytic function. The
function $h_0$ can be extended to, generally speaking multi-valued,
analytic function ${\cal A}$  along any path in $D_*$ because $u$ is
given in the whole domain  $D_*$. Thus, $f={\cal A}\circ g$ is a
desired multi-valued solution of the Dirichlet problem
(\ref{eqDIRICHLET}) for Beltrami equation (\ref{eqBeltrami}).
\end{proof} $\Box$

\medskip

{\bf Remark 9.} Note that if the family of the functions
$\psi_{z_0,\varepsilon}(t)\equiv\psi_{z_0}(t)$ is independent on the
parameter $\varepsilon$, then the condition (\ref{3omalM}) implies
that $I_{z_0}(\varepsilon)\to \infty$ as $\varepsilon\to 0$. This
follows immediately from arguments by contradiction, apply for it
(\ref{eqConnect}) and the condition $K_{\mu}\in L^1(D)$. Note also
that (\ref{3omalM}) holds, in particular, if, for some
$\varepsilon_0=\varepsilon(z_0)$,
\begin{equation}\label{333omalM}
\int\limits_{|z-z_0|<\varepsilon_0}
K^T_{\mu}(z,z_0)\cdot\psi_{z_0}^2(|z-z_0|)\,dm(z)<\infty \qquad
\forall\ z_0\in\overline{D}\end{equation} and
$I_{z_0}(\varepsilon)\to \infty$ as $\varepsilon\to 0$. In other
words, for the existence of a multi-valued solutions for the
Dirichlet problem (\ref{eqDIRICHLET}) to the Beltrami equation
(\ref{eqBeltrami}) with each continuous boundary functions
$\varphi$, it is sufficient that the integral in (\ref{333omalM})
converges for some nonnegative function $\psi_{z_0}(t)$ that is
locally integrable over $(0,\varepsilon_0 ]$ but has a nonintegrable
singularity at $0$. The functions $\log^{\lambda}(e/|z-z_0|)$,
$\lambda\in (0,1)$, $z\in\Bbb D$, $z_0\in\overline{\Bbb D}$, and
$\psi(t)=1/(t \,\, \log(e/t))$, $t\in(0,1)$, show that the condition
(\ref{333omalM}) is compatible with the condition
$I_{z_0}(\varepsilon)\to\infty $ as $\varepsilon\to 0$. Furthermore,
the condition (\ref{3omalM}) in Lemma 5 shows that it is sufficient
for the existence of a multi-valued solutions for the Dirichlet
problem (\ref{eqDIRICHLET}) to the Beltrami equation
(\ref{eqBeltrami}) even that the integral in (\ref{333omalM}) to be
divergent in a controlled way.

\medskip

Arguing similarly to the last section, we derive from Lemma 5 the
new series of results.

\medskip

For instance, choosing $\psi(t)=1/\left(t\,
\log\left(1/t\right)\right)$ in Lemma 5, we obtain by Lemma 2 the
following.

\medskip

{\bf Theorem 5.} {\it Let $D$ be a bounded domain in ${\Bbb C}$ with
no boundary component degenerated to a single point and
$\mu:D\to{\Bbb C}$ be a measurable function with $|\mu(z)|<1$ a.e.
and $K_{\mu}\in L^1(D)$. Suppose that $K^T_{\mu}(z,z_0)\leqslant
Q_{z_0}(z)$ a.e. in $U_{z_0}$ for every point $z_0\in \overline{D}$,
a neighborhood $U_{z_0}$ of $z_0$ and a function $Q_{z_0}:
U_{z_0}\to[0,\infty]$ in the class ${\rm FMO}({z_0})$. Then the
Beltrami equation (\ref{eqBeltrami}) has a multi-valued solution of
the Dirichlet problem (\ref{eqDIRICHLET}) with the representation
(\ref{eqHYDROm}) for each continuous function $\varphi:\partial
D\to{\Bbb R}$.}

\medskip

In particular, by Proposition 1 the conclusion of Theorem 5 holds if
every point $z_0\in\overline{D}$ is the Lebesgue point of a suitable
dominant $Q_{z_0}$.

\bigskip

By Corollary 2 we obtain the following nice consequence of Theorem
5, too.

\medskip

{\bf Corollary 11.} {\it  Let $D$ be a bounded domain in ${\Bbb C}$
with no boundary component degenerated to a single point and
$\mu:D\to{\Bbb C}$ be a measurable function with $|\mu(z)|<1$ a.e.,
$K_{\mu}\in L^1(D)$ and
\begin{equation}\label{eqMEANm}\overline{\lim\limits_{\varepsilon\to0}}\quad
\dashint_{B(z_0,\varepsilon)}K^T_{\mu}(z,z_0)\,dm(z)<\infty\qquad\forall\
z_0\in\overline{D}\, .\end{equation}Then the Beltrami equation
(\ref{eqBeltrami}) has a multi-valued solution of the Dirichlet
problem (\ref{eqDIRICHLET}) with the representation (\ref{eqHYDROm})
for each continuous function $\varphi:\partial D\to{\Bbb R}$.}

\medskip

Since $K^T_{\mu}(z,z_0) \leqslant K_{\mu}(z)$ for all $z$ and
$z_0\in \Bbb C$, we also obtain the following consequences of
Theorem 5.

\medskip

{\bf Corollary 12.} {\it Let $D$ be a bounded domain in ${\Bbb C}$
with no boundary component degenerated to a single point and
$\mu:D\to{\Bbb C}$ be a measurable function with $|\mu(z)|<1$ a.e.
and $K_{\mu}$ have a dominant $Q:\mathbb C\to[1,\infty)$ in the
class {\rm BMO}$_{\rm loc}$. Then the Beltrami equation
(\ref{eqBeltrami}) has a multi-valued solution of the Dirichlet
problem (\ref{eqDIRICHLET}) with the representation (\ref{eqHYDROm})
for each continuous function $\varphi:\partial D\to{\Bbb R}$.}

\medskip

{\bf Remark 10.} In particular, the conclusion of Corollary 12 holds
if $Q\in{\rm W}^{1,2}_{\rm loc}$ because $W^{\,1,2}_{\rm loc}
\subset {\rm VMO}_{\rm loc}$, see \cite{BN}.

\medskip

{\bf Corollary 13.} {\it Let $D$ be a bounded domain in ${\Bbb C}$
with no boundary component degenerated to a single point and
$\mu:D\to{\Bbb C}$ be a measurable function with $|\mu(z)|<1$ and
$K_{\mu}(z)\leqslant Q(z)$ a.e. in $D$ with a function $Q$ in the
class ${\rm FMO}(\overline{D})$. Then the Beltrami equation
(\ref{eqBeltrami}) has a multi-valued solution of the Dirichlet
problem (\ref{eqDIRICHLET}) with the representation (\ref{eqHYDROm})
for each continuous function $\varphi:\partial D\to{\Bbb R}$.}

\medskip

Similarly, choosing in Lemma 5 the function $\psi(t)=1/t$, we come
to the next statement.

\medskip

{\bf Theorem 6.} {\it Let $D$ be a bounded domain in ${\Bbb C}$ with
no boundary component degenerated to a single point and
$\mu:D\to{\Bbb C}$ be a measurable function with $|\mu(z)|<1$ a.e.
and $K_{\mu}\in L^1(D)$. Suppose that
\begin{equation}\label{eqLOGm}
\int\limits_{\varepsilon<|z-z_0|<\varepsilon_0}K^T_{\mu}(z,z_0)\,\frac{dm(z)}{|z-z_0|^2}
=o\left(\left[\log\frac{1}{\varepsilon}\right]^2\right)\qquad\hbox{as
$\varepsilon\to 0$}\qquad\forall\ z_0\in\overline{D}\end{equation}
for some $\varepsilon_0=\varepsilon(z_0)>0$. Then Beltrami equation
(\ref{eqBeltrami}) has a multi-valued solution of the Dirichlet
problem (\ref{eqDIRICHLET}) with representation (\ref{eqHYDROm}) for
each continuous function $\varphi:\partial D\to{\Bbb R}$.}

\medskip

{\bf Remark 11.} Choosing in Lemma 5 the function
$\psi(t)=1/(t\log{1/t})$ instead of $\psi(t)=1/t$, we are able to
replace (\ref{eqLOGm}) by
\begin{equation}\label{eqLOGLOGm}
\int\limits_{\varepsilon<|z-z_0|<\varepsilon_0}\frac{K^T_{\mu}(z,z_0)\,dm(z)}
{\left(|z-z_0|\log{\frac{1}{|z-z_0|}}\right)^2}
=o\left(\left[\log\log\frac{1}{\varepsilon}\right]^2\right)\end{equation}
In general, we are able to give here the whole scale of the
corresponding conditions in $\log$ using functions $\psi(t)$ of the
form
$1/(t\log{1}/{t}\cdot\log\log{1}/{t}\cdot\ldots\cdot\log\ldots\log{1}/{t})$.

\medskip

Choosing in Lemma 5 the functional parameter
${\psi}_{z_0,{\varepsilon}}(t) \equiv {\psi}_{z_0}(t) \colon =
1/[tk^T_{\mu}(z_0,t)]$, where $k_{\mu}^T(z_0,r)$ is the integral
mean of $K^T_{{\mu}}(z,z_0)$ over the circle $S(z_0,r)\, :=\, \{ z
\in\mathbb C:\, |z-z_0|\,=\, r\}$, we obtain one more important
conclusion.

\medskip

{\bf Theorem 7.} {\it Let $D$ be a bounded domain in ${\Bbb C}$ with
no boundary component degenerated to a single point and
$\mu:D\to{\Bbb C}$ be a measurable function with $|\mu(z)|<1$ a.e.
and $K_{\mu}\in L^1(D)$. Suppose that
\begin{equation}\label{eqLEHTOm}\int\limits_{0}^{\varepsilon_0}
\frac{dr}{rk^T_{\mu}(z_0,r)}=\infty\qquad\forall\
z_0\in\overline{D}\end{equation} for some
$\varepsilon_0=\varepsilon(z_0)>0$. Then Beltrami equation
(\ref{eqBeltrami}) has a multi-valued solution of the Dirichlet
problem (\ref{eqDIRICHLET}) with representation (\ref{eqHYDROm}) for
each continuous function $\varphi:\partial D\to{\Bbb R}$.}

\medskip

{\bf Corollary 14.} {\it Let $D$ be a bounded domain in ${\Bbb C}$
with no boundary component degenerated to a single point and
$\mu:D\to{\Bbb C}$ be a measurable function with $|\mu(z)|<1$ a.e.,
$K_{\mu}\in L^1(D)$ and
\begin{equation}\label{eqLOGkM}k^T_{\mu}(z_0,\varepsilon)=O\left(\log\frac{1}{\varepsilon}\right)
\qquad\mbox{as}\ \varepsilon\to0\qquad\forall\ z_0\in\overline{D}\
.\end{equation} Then the Beltrami equation (\ref{eqBeltrami}) has a
multi-valued solution of the Dirichlet problem (\ref{eqDIRICHLET})
with the representation (\ref{eqHYDROm}) for each continuous
function $\varphi:\partial D\to{\Bbb R}$.}

\medskip

{\bf Remark 12.} In particular, the conclusion of Corollary 14 holds
if
\begin{equation}\label{eqLOGKm} K^T_{\mu}(z,z_0)=O\left(\log\frac{1}{|z-z_0|}\right)\qquad{\rm
as}\quad z\to z_0\quad\forall\ z_0\in\overline{D}\,.\end{equation}
Moreover, the condition (\ref{eqLOGkM}) can be replaced by the whole
series of more weak conditions
\begin{equation}\label{edLOGLOGkM}
k^T_{\mu}(z_0,\varepsilon)=O\left(\left[\log\frac{1}{\varepsilon}\cdot\log\log\frac{1}
{\varepsilon}\cdot\ldots\cdot\log\ldots\log\frac{1}{\varepsilon}
\right]\right) \qquad\forall\ z_0\in \overline{D}\ .
\end{equation}

\medskip

Combining Theorems 7, Proposition 2 and Remark 3, we obtain the
following result.

\medskip

{\bf Theorem 8.} {\it Let $D$ be a bounded domain in ${\Bbb C}$ with
no boundary component degenerated to a single point and
$\mu:D\to{\Bbb C}$ be a measurable function with $|\mu(z)|<1$ a.e.
and $K_{\mu}\in L^1(D)$. Suppose that
\begin{equation}\label{eqINTEGRALm}\int\limits_{U_{z_0}}\Phi_{z_0}\left(K^T_{\mu}(z,z_0)\right)\,dm(z)<\infty
\qquad\forall\ z_0\in \overline{D}\end{equation} for a neighborhood
$U_{z_0}$ of $z_0$ and a convex non-decreasing function
$\Phi_{z_0}:[0,\infty]\to[0,\infty]$ with
\begin{equation}\label{eqINTm}
\int\limits_{\Delta(z_0)}^{\infty}\log\,\Phi_{z_0}(t)\,\frac{dt}{t^2}\
=\ +\infty\end{equation} for some $\Delta(z_0)>0$. Then Beltrami
equation (\ref{eqBeltrami}) has a multi-valued solution of the
Dirichlet problem (\ref{eqDIRICHLET}) with the representation
(\ref{eqHYDROm}) for each continuous function $\varphi:\partial
D\to{\Bbb R}$.}

\medskip

{\bf Corollary 15.} {\it Let $D$ be a bounded domain in ${\Bbb C}$
with no boundary component degenerated to a single point and
$\mu:D\to{\Bbb C}$ be a measurable function with $|\mu(z)|<1$ a.e.,
$K_{\mu}\in L^1(D)$ and
\begin{equation}\label{eqEXPm}\int\limits_{U_{z_0}}e^{\alpha(z_0) K^T_{\mu}(z,z_0)}\,dm(z)<\infty
\qquad\forall\ z_0\in \overline{D}\end{equation} for some
$\alpha(z_0)>0$ and a neighborhood $U_{z_0}$ of the point $z_0$.
Then the Beltrami equation (\ref{eqBeltrami}) has a multi-valued
solution of the Dirichlet problem (\ref{eqDIRICHLET}) with the
representation (\ref{eqHYDROm}) for each continuous function
$\varphi:\partial D\to{\Bbb R}$.}

\medskip

Since $K^T_{\mu}(z,z_0) \leqslant K_{\mu}(z)$ for $z$ and $z_0\in
\Bbb C$ and $z\in D$, we also obtain the following consequences of
Theorem 8.

\medskip

{\bf Corollary 16.} {\it Let $D$ be a bounded domain in ${\Bbb C}$
with no boundary component degenerated to a single point and
$\mu:D\to{\Bbb C}$ be a measurable function with $|\mu(z)|<1$ a.e.
and $K_{\mu}\in L^1(D)$. Suppose that
\begin{equation}\label{eqINTKm}\int\limits_{D}\Phi\left(K_{\mu}(z)\right)\,dm(z)<\infty\end{equation}
for a convex non-decreasing function $\Phi:[0,\infty]\to[0,\infty]$
with
\begin{equation}\label{eqINTFm}
\int\limits_{\delta}^{\infty}\log\,\Phi(t)\,\frac{dt}{t^2}\ =\
+\infty\end{equation} for some $\delta>0$. Then Beltrami equation
(\ref{eqBeltrami}) has a multi-valued solution of the Dirichlet
problem (\ref{eqDIRICHLET}) with the representation (\ref{eqHYDROm})
for each continuous function $\varphi:\partial D\to{\Bbb R}$.}

\medskip

{\bf Corollary 17.} {\it Let $D$ be a bounded domain in ${\Bbb C}$
with no boundary component degenerated to a single point and
$\mu:D\to{\Bbb C}$ be a measurable function with $|\mu(z)|<1$ a.e.
and, for some $\alpha>0$,
\begin{equation}\label{eqEXPAm}\int\limits_{D}e^{\alpha K_{\mu}(z)}\,dm(z)\ <\
\infty\ .
\end{equation} Then the Beltrami equation
(\ref{eqBeltrami}) has a multi-valued solution of the Dirichlet
problem (\ref{eqDIRICHLET}) with the representation (\ref{eqHYDROm})
for each continuous function $\varphi:\partial D\to{\Bbb R}$.}

\medskip

{\bf Remark 13.} By the Stoilow theorem, see e.g. \cite{Sto}, a
multi-valued solution $f$ of the Dirichlet problem
(\ref{eqDIRICHLET}) for the Beltrami equation (\ref{eqBeltrami})
with $K_{\mu}\in L^1_{\rm loc}(D)$ can be represented in the form
$f={\cal A}\circ F$ where $\cal A$ is a multi-valued analytic
function and $F$ is a homeomorphic regular solution of
(\ref{eqBeltrami}) in the class $W_{\rm loc}^{1,1}$. Thus, by
Theorem 5.1 in \cite{RSY$_5$} the condition (\ref{eqINTFm}) is not
only sufficient but also necessary to have a regular solution of the
Dirichlet problem (\ref{eqDIRICHLET}) for arbitrary Beltrami
equations (\ref{eqBeltrami}) with the integral constraints
(\ref{eqINTKm}) for all continuous functions $\varphi:\partial
D\to\Bbb{R}$, see also Remark 3.

\section{Applications to the potential theory}

The results of the last section seem too abstract, and therefore
allegedly useless. However, we give here their significant
applications to one of the main equations of mathematical physics in
anisotropic and inhomogeneous media.

Namely, in this section we obtain criteria for the existence and
representation of solutions $u$ of the Dirichlet problem to the
elliptic equations of the form
\begin{equation}\label{eqPotential}
{\rm div}\, A\nabla\,u=0
\end{equation}
with measurable matrix-valued function $A(z)=\{a_{ij}(z)\}$ in
arbitrary bounded domains $D$ in $\mathbb C$ with no boundary
component degenerated to a single point. Our example in Introduction
shows that such domains form the most wide class of domains for
which the Dirichlet problem will be always solvable for each
continuous boundary function $\varphi :\partial D\to\mathbb R$.

Solutions of the equation (\ref{eqPotential}) are called {\bf
$A-$harmonic functions}, see e.g. \cite{HKM}, and they satisfied
(\ref{eqPotential}) in the sense of distributions, i.e., in the
sense that $u\in W^{1,1}_{\rm loc}(D)$ and that
\begin{equation}\label{weak}
\int_D \langle A(z)\nabla u(z),\nabla\psi(z)\rangle\
d\,m(z)=0\,\,\,\,\,\,\,\,\,\forall\,\,\, \psi\in C_0^\infty(D)\ ,
\end{equation}
where $C_0^\infty(D)$ denotes the collection of all infinitely
differentiable functions $\psi: D\to\mathbb R$ with compact support
in $D$, $\langle a,b\rangle$ means the scalar product of vectors $a$
and $b$ in $\mathbb R^2$, and $d\,m(z):=d\,x\,d\,y$, $z=x+iy$,
corresponds to the Lebesgue measure (area) in the plane $\mathbb C$.

In this connection, let us describe the relevance of the Beltrami
equations (\ref{eqBeltrami}) and the equations (\ref{eqPotential}).
First of all, recall that the {\bf Hodge operator} $J$ is the
counterclockwise rotation by the angle $\pi/2$ in $\mathbb R^2$:
\begin{equation}\label{eqHodge}
J=\left[\begin{array}{rr} 0 &-1 \\
        1 &0\end{array}\right]\ :\ \mathbb R^2\to\mathbb R^2\ ,\ \ \ \
        \ \ \ J^2\ =\ -\ I\ ,
\end{equation}
where $I$ denotes the unit $2\times 2$ matrix. Thus, the matrix $J$
plays the role of an imaginary unit in the space of two-dimensional
square matrices with real-valued elements.

By Theorem 16.1.6 in \cite{AIM}, if $f$ is a $W^{1,1}_{\rm loc}$
solution of the Beltrami equation (\ref{eqBeltrami}), then the
functions $u:={\rm Re} f$ and $v:={\rm Im} f$ satisfy the equation:
\begin{equation}\label{eqAconjugate}
\nabla\, v(z)\ =\ J\, A(z)\nabla\,u(z)\ ,
\end{equation}
where the matrix-valued function $A(z)$ is calculated through
$\mu(z)$ in the following way:
\begin{equation}\label{matrix}
A\  =\ \left[\begin{array}{ccc} a_{11}  & a_{12} \\
a_{21} & a_{22} \end{array}\right]\ :=\
\left[\begin{array}{ccc} {|1-\mu|^2\over 1-|\mu|^2}  & {-2{\rm Im}\,\mu\over 1-|\mu|^2} \\
{-2{\rm Im}\,\mu\over 1-|\mu|^2} & {|1+\mu|^2\over 1-|\mu|^2}
\end{array}\right].\end{equation}
The function $v$ is called the {\bf $A-$harmonic conjugate of} u or
sometimes a {\bf stream function of the potential} $u$. Note that by
(\ref{eqHodge}) the equation (\ref{eqAconjugate}) is equivalent to
the equation
\begin{equation}\label{eqDIV}
A(z)\nabla\,u(z)\ =\ -J\nabla\, v(z)\ .
\end{equation}
As known, the curl of any gradient field is zero in the sense of
distributions and the Hodge operator $J$ transforms curl-free fields
into divergence-free fields, and vice versa, see e.g. 16.1.3 in
\cite{AIM}. Hence (\ref{eqDIV}) implies (\ref{eqPotential}).

We see from (\ref{matrix}) that the matrix $A$ is symmetric and it
is clear by elementary calculations that ${\rm det}\,A=1$. Moreover,
since $|\mu(z)|<1$ a.e., from ellipticity of this matrix $A$ follows
that ${\rm det}\,(I+A)>0$ a.e., which in terms of its elements means
that $(1+a_{11})(1+a_{22})>a_{12}a_{21}$ a.e. Further $\mathbb
S^{2\times 2}$ denotes the collection of all such matrices. Thus, by
Theorem 16.1.6 in \cite{AIM}, the Beltrami equation is the complex
form of one of the main equations of mathematical physics, the
potential equation (\ref{eqPotential}) with the matrix-valued
coefficient $A$ in the class $\mathbb S^{2\times 2}$.

Note that the matrix identities in (\ref{matrix}) can be converted
a.e. to express the coefficient $\mu(z)$ of the Beltrami equation
(\ref{eqBeltrami}) through the elements of the matrices $A(z)$:
\begin{equation}\label{eqCoefficient}
\mu\ =\ \mu_A\ :=\
-\frac{a_{11}-a_{22}+i(a_{12}+a_{21})}{2+a_{11}+a_{22}}\ .
\end{equation}
Thus, we obtain the latter expression as a criterion for the
solvability of the Dirichlet problem
\begin{equation}\label{eqDIR}
\lim_{z\to\zeta}\ u(z)\ =\ \varphi(\zeta)\ \ \ \ \ \forall\
\zeta\in\partial D
\end{equation}
to the potential equation (\ref{eqPotential}). Namely, by the above
arguments in this section as well as Lemma 5 we come to the
following general criteria.

\medskip

{\bf Lemma 6.} {\it Let $D$ be a bounded domain in ${\Bbb C}$ with
no boundary component degenerated to a single point, $A:D\to\mathbb
S^{2\times 2}$ be a measurable function in $D$ with $K_{\mu_A}\in
L^1(D)$ and
\begin{equation}\label{3omalMA}
\int\limits_{\varepsilon<|z-z_0|<\varepsilon_0}
K^T_{\mu_A}(z,z_0)\cdot\psi^2_{z_0,\varepsilon}(|z-z_0|)\,dm(z)=o(I_{z_0}^{2}(\varepsilon))\quad{\rm
as}\quad\varepsilon\to0\ \ \forall\ z_0\in\overline{D}\end{equation}
for some $\varepsilon_0=\varepsilon(z_0)>0$ and a family of
measurable functions $\psi_{z_0,\varepsilon}:
(0,\varepsilon_0)\to(0,\infty)$ with
\begin{equation}\label{eq3.5.3MA}
I_{z_0}(\varepsilon)\colon
=\int\limits_{\varepsilon}^{\varepsilon_0}
\psi_{z_0,\varepsilon}(t)\,dt<\infty\qquad\forall\
\varepsilon\in(0,\varepsilon_0)\,.\end{equation} Then the potential
equation (\ref{eqPotential}) has $A-$harmonic solutions $u$ of the
Dirichlet problem (\ref{eqDIR}) for each continuous function
$\varphi:\partial D\to{\Bbb R}$.

Moreover, such a solution $u$ can be represented as the composition
\begin{equation}\label{eqHYDROmA}
u\ =\ {\cal H}\circ g\ ,\ \ \ \ \ \ \hbox{$g(z)\, =\ z\, +\, o(1)$\
\ \ as\ \ \ $z\to\infty$}\ ,
\end{equation}
where $g:\mathbb C\to\mathbb C$ is a regular homeomorphic solution
of the Beltrami equation (\ref{eqBeltrami}) in $\mathbb C$ with
$\mu_A$ extended by zero outside of $D$ and ${\cal H}:D_*\to\mathbb
C,$ $D_*:=g(D)$, is a unique harmonic function satisfying the
Dirichlet condition
\begin{equation}\label{eqHOLOMORPHICMA}
\lim_{\xi\to\zeta}\ {\cal H}(\xi)\ =\ \varphi_*(\zeta)\ \ \ \ \
\forall\ \zeta\in\partial D_*\ ,\ \ \ \ \ \mbox{where
$\varphi_*:=\varphi\circ g^{-1}$.}
\end{equation}}

As it was before, we assume here that the dilatations
$K^T_{\mu_A}(z,z_0)$ and $K_{\mu_A}(z)$ are extended by $1$ outside
of the domain $D$.

\medskip

{\bf Remark 14.} Note that if the family of the functions
$\psi_{z_0,\varepsilon}(t)\equiv\psi_{z_0}(t)$ is independent on the
parameter $\varepsilon$, then the condition (\ref{3omalMA}) implies
that $I_{z_0}(\varepsilon)\to \infty$ as $\varepsilon\to 0$. This
follows immediately from arguments by contradiction, apply for it
(\ref{eqConnect}) and the condition $K_{\mu_A}\in L^1(D)$. Note also
that (\ref{3omalMA}) holds, in particular, if, for some
$\varepsilon_0=\varepsilon(z_0)$,
\begin{equation}\label{333omalMA}
\int\limits_{|z-z_0|<\varepsilon_0}
K^T_{\mu_A}(z,z_0)\cdot\psi_{z_0}^2(|z-z_0|)\,dm(z)<\infty \qquad
\forall\ z_0\in\overline{D}\end{equation} and
$I_{z_0}(\varepsilon)\to \infty$ as $\varepsilon\to 0$. In other
words, for the existence of $A-$harmonic solutions of the Dirichlet
problem (\ref{eqDIR}) to the potential equation (\ref{eqPotential})
with each continuous boundary functions $\varphi$, it is sufficient
that the integral in (\ref{333omalMA}) converges for some
nonnegative function $\psi_{z_0}(t)$ that is locally integrable over
$(0,\varepsilon_0 ]$ but has a nonintegrable singularity at $0$. The
functions $\log^{\lambda}(e/|z-z_0|)$, $\lambda\in (0,1)$, $z\in\Bbb
D$, $z_0\in\overline{\Bbb D}$, and $\psi(t)=1/(t \,\, \log(e/t))$,
$t\in(0,1)$, show that the condition (\ref{333omalMA}) is compatible
with the condition $I_{z_0}(\varepsilon)\to\infty $ as
$\varepsilon\to 0$. Furthermore, the condition (\ref{3omalMA}) in
Lemma 6 shows that it is sufficient for the existence of
$A-$harmonic solutions of the Dirichlet problem (\ref{eqDIR}) to the
potential equation (\ref{eqPotential}) even that the integral in
(\ref{333omalMA}) to be divergent in a controlled way.

\medskip

Arguing similarly to Section 4, we derive from Lemma 6 the next
series of results.

\medskip

For instance, choosing $\psi(t)=1/\left(t\,
\log\left(1/t\right)\right)$ in Lemma 6, we obtain by Lemma 2 the
following.

\medskip

{\bf Theorem 9.} {\it Let $D$ be a bounded domain in ${\Bbb C}$ with
no boundary component degenerated to a single point and
$A:D\to\mathbb S^{2\times 2}$ be a measurable function in $D$ with
$K_{\mu_A}\in L^1(D)$. Suppose that $K^T_{\mu_A}(z,z_0)\leqslant
Q_{z_0}(z)$ a.e. in $U_{z_0}$ for every point $z_0\in \overline{D}$,
a neighborhood $U_{z_0}$ of $z_0$ and a function $Q_{z_0}:
U_{z_0}\to[0,\infty]$ in the class ${\rm FMO}({z_0})$. Then the
potential equation (\ref{eqPotential}) has $A-$harmonic solutions of
the Dirichlet problem (\ref{eqDIR}) with the representation
(\ref{eqHYDROmA}) for each continuous function $\varphi:\partial
D\to{\Bbb R}$.}

\medskip

In particular, by Proposition 1 the conclusion of Theorem 9 holds if
every point $z_0\in\overline{D}$ is the Lebesgue point of a suitable
dominant $Q_{z_0}$.

\bigskip

By Corollary 2 we obtain the following nice consequence of Theorem
9, too.

\medskip

{\bf Corollary 18.} {\it  Let $D$ be a bounded domain in ${\Bbb C}$
with no boundary component degenerated to a single point and
$A:D\to\mathbb S^{2\times 2}$ be a measurable function in $D$ with
$K_{\mu_A}\in L^1(D)$ and
\begin{equation}\label{eqMEANmA}\overline{\lim\limits_{\varepsilon\to0}}\quad
\dashint_{B(z_0,\varepsilon)}K^T_{\mu_A}(z,z_0)\,dm(z)<\infty\qquad\forall\
z_0\in\overline{D}\, .\end{equation} Then the potential equation
(\ref{eqPotential}) has $A-$harmonic solutions of the Dirichlet
problem (\ref{eqDIR}) with the representation (\ref{eqHYDROmA}) for
each continuous function $\varphi:\partial D\to{\Bbb R}$.}

\medskip

Since $K^T_{\mu_A}(z,z_0) \leqslant K_{\mu_A}(z)$ for all $z$ and
$z_0\in \Bbb C$, we also obtain the following consequences of
Theorem 9.

\medskip

{\bf Corollary 19.} {\it Let $D$ be a bounded domain in ${\Bbb C}$
with no boundary component degenerated to a single point,
$A:D\to\mathbb S^{2\times 2}$ be a measurable function in $D$ and
$K_{\mu_A}$ have a dominant $Q:\mathbb C\to[1,\infty)$ in the class
{\rm BMO}$_{\rm loc}$. Then the potential equation
(\ref{eqPotential}) has $A-$harmonic solutions of the Dirichlet
problem (\ref{eqDIR}) with the representation (\ref{eqHYDROmA}) for
each continuous function $\varphi:\partial D\to{\Bbb R}$.}

\medskip

{\bf Remark 15.} In particular, the conclusion of Corollary 19 holds
if $Q\in{\rm W}^{1,2}_{\rm loc}$ because $W^{\,1,2}_{\rm loc}
\subset {\rm VMO}_{\rm loc}$, see \cite{BN}.

\medskip

{\bf Corollary 20.} {\it Let $D$ be a bounded domain in ${\Bbb C}$
with no boundary component degenerated to a single point,
$A:D\to\mathbb S^{2\times 2}$ be a measurable function in $D$ and
$K_{\mu_A}(z)\leqslant Q(z)$ a.e. in $D$ with a function $Q$ in the
class ${\rm FMO}(\overline{D})$. Then the potential equation
(\ref{eqPotential}) has $A-$harmonic solutions of the Dirichlet
problem (\ref{eqDIR}) with the representation (\ref{eqHYDROmA}) for
each continuous function $\varphi:\partial D\to{\Bbb R}$.}

\medskip

Similarly, choosing in Lemma 6 the function $\psi(t)=1/t$, we come
to the next statement.

\medskip

{\bf Theorem 10.} {\it Let $D$ be a bounded domain in ${\Bbb C}$
with no boundary component degenerated to a single point,
$A:D\to\mathbb S^{2\times 2}$ be a measurable function in $D$ with
$K_{\mu_A}\in L^1(D)$. Suppose that
\begin{equation}\label{eqLOGmA}
\int\limits_{\varepsilon<|z-z_0|<\varepsilon_0}K^T_{\mu_A}(z,z_0)\,\frac{dm(z)}{|z-z_0|^2}
=o\left(\left[\log\frac{1}{\varepsilon}\right]^2\right)\qquad\hbox{as
$\varepsilon\to 0$}\qquad\forall\ z_0\in\overline{D}\end{equation}
for some $\varepsilon_0=\varepsilon(z_0)>0$. Then the potential
equation (\ref{eqPotential}) has $A-$harmonic solutions of the
Dirichlet problem (\ref{eqDIR}) with representation
(\ref{eqHYDROmA}) for each continuous function $\varphi:\partial
D\to{\Bbb R}$.}

\medskip

{\bf Remark 16.} Choosing in Lemma 6 the function
$\psi(t)=1/(t\log{1/t})$ instead of $\psi(t)=1/t$, we are able to
replace (\ref{eqLOGmA}) by
\begin{equation}\label{eqLOGLOGmA}
\int\limits_{\varepsilon<|z-z_0|<\varepsilon_0}\frac{K^T_{\mu_A}(z,z_0)\,dm(z)}
{\left(|z-z_0|\log{\frac{1}{|z-z_0|}}\right)^2}
=o\left(\left[\log\log\frac{1}{\varepsilon}\right]^2\right)\end{equation}
In general, we are able to give here the whole scale of the
corresponding conditions in $\log$ using functions $\psi(t)$ of the
form
$1/(t\log{1}/{t}\cdot\log\log{1}/{t}\cdot\ldots\cdot\log\ldots\log{1}/{t})$.

\medskip

Choosing in Lemma 6 the functional parameter
${\psi}_{z_0,{\varepsilon}}(t) \equiv {\psi}_{z_0}(t) \colon =
1/[tk^T_{\mu_A}(z_0,t)]$, where $k_{\mu_A}^T(z_0,r)$ is the integral
mean of $K^T_{{\mu_A}}(z,z_0)$ over the circle $S(z_0,r)\, :=\, \{ z
\in\mathbb C:\, |z-z_0|\,=\, r\}$, we obtain one more important
conclusion.

\medskip

{\bf Theorem 11.} {\it Let $D$ be a bounded domain in ${\Bbb C}$
with no boundary component degenerated to a single point,
$A:D\to\mathbb S^{2\times 2}$ be a measurable function in $D$ with
$K_{\mu_A}\in L^1(D)$. Suppose that
\begin{equation}\label{eqLEHTOmA}\int\limits_{0}^{\varepsilon_0}
\frac{dr}{rk^T_{\mu_A}(z_0,r)}=\infty\qquad\forall\
z_0\in\overline{D}\end{equation} for some
$\varepsilon_0=\varepsilon(z_0)>0$. Then the potential equation
(\ref{eqPotential}) has $A-$harmonic solutions of the Dirichlet
problem (\ref{eqDIR}) with representation (\ref{eqHYDROmA}) for each
continuous function $\varphi:\partial D\to{\Bbb R}$.}

\medskip

{\bf Corollary 21.} {\it Let $D$ be a bounded domain in ${\Bbb C}$
with no boundary component degenerated to a single point,
$A:D\to\mathbb S^{2\times 2}$ be a measurable function in $D$ with
$K_{\mu_A}\in L^1(D)$ and
\begin{equation}\label{eqLOGkMA}k^T_{\mu_A}(z_0,\varepsilon)=O\left(\log\frac{1}{\varepsilon}\right)
\qquad\mbox{as}\ \varepsilon\to0\qquad\forall\ z_0\in\overline{D}\
.\end{equation} Then the potential equation (\ref{eqPotential}) has
$A-$harmonic solutions of the Dirichlet problem (\ref{eqDIR}) with
the representation (\ref{eqHYDROmA}) for each continuous function
$\varphi:\partial D\to{\Bbb R}$.}

\medskip

{\bf Remark 17.} In particular, the conclusion of Corollary 21 holds
if
\begin{equation}\label{eqLOGKmA} K^T_{\mu_A}(z,z_0)=O\left(\log\frac{1}{|z-z_0|}\right)\qquad{\rm
as}\quad z\to z_0\quad\forall\ z_0\in\overline{D}\,.\end{equation}
Moreover, the condition (\ref{eqLOGkMA}) can be replaced by the
whole series of more weak conditions
\begin{equation}\label{edLOGLOGkMA}
k^T_{\mu_A}(z_0,\varepsilon)=O\left(\left[\log\frac{1}{\varepsilon}\cdot\log\log\frac{1}
{\varepsilon}\cdot\ldots\cdot\log\ldots\log\frac{1}{\varepsilon}
\right]\right) \qquad\forall\ z_0\in \overline{D}\ .
\end{equation}

\medskip

Combining Theorems 11, Proposition 2 and Remark 3, we obtain the
following result.

\medskip

{\bf Theorem 12.} {\it Let $D$ be a bounded domain in ${\Bbb C}$
with no boundary component degenerated to a single point,
$A:D\to\mathbb S^{2\times 2}$ be a measurable function in $D$ with
$K_{\mu_A}\in L^1(D)$. Suppose that
\begin{equation}\label{eqINTEGRALmA}\int\limits_{U_{z_0}}\Phi_{z_0}\left(K^T_{\mu_A}(z,z_0)\right)\,dm(z)<\infty
\qquad\forall\ z_0\in \overline{D}\end{equation} for a neighborhood
$U_{z_0}$ of $z_0$ and a convex non-decreasing function
$\Phi_{z_0}:[0,\infty]\to[0,\infty]$ with
\begin{equation}\label{eqINTmA}
\int\limits_{\Delta(z_0)}^{\infty}\log\,\Phi_{z_0}(t)\,\frac{dt}{t^2}\
=\ +\infty\end{equation} for some $\Delta(z_0)>0$. Then the
potential equation (\ref{eqPotential}) has $A-$harmonic solutions of
the Dirichlet problem (\ref{eqDIR}) with representation
(\ref{eqHYDROmA}) for each continuous function $\varphi:\partial
D\to{\Bbb R}$.}

\medskip

{\bf Corollary 22.} {\it Let $D$ be a bounded domain in ${\Bbb C}$
with no boundary component degenerated to a single point,
$A:D\to\mathbb S^{2\times 2}$ be a measurable function in $D$ with
$K_{\mu_A}\in L^1(D)$ and
\begin{equation}\label{eqEXPmA}\int\limits_{U_{z_0}}e^{\alpha(z_0) K^T_{\mu_A}(z,z_0)}\,dm(z)<\infty
\qquad\forall\ z_0\in \overline{D}\end{equation} for some
$\alpha(z_0)>0$ and a neighborhood $U_{z_0}$ of the point $z_0$.
Then the potential equation (\ref{eqPotential}) has $A-$harmonic
solutions of the Dirichlet problem (\ref{eqDIR}) with the
representation (\ref{eqHYDROmA}) for each continuous function
$\varphi:\partial D\to{\Bbb R}$.}

\medskip

Since $K^T_{\mu_A}(z,z_0) \leqslant K_{\mu_A}(z)$ for $z$ and
$z_0\in \Bbb C$ and $z\in D$, we also obtain the following
consequences of Theorem 12.

\medskip

{\bf Corollary 23.} {\it Let $D$ be a bounded domain in ${\Bbb C}$
with no boundary component degenerated to a single point,
$A:D\to\mathbb S^{2\times 2}$ be a measurable function in $D$ with
$K_{\mu_A}\in L^1(D)$. Suppose that
\begin{equation}\label{eqINTKmA}\int\limits_{D}\Phi\left(K_{\mu_A}(z)\right)\,dm(z)<\infty\end{equation}
for a convex non-decreasing function $\Phi:[0,\infty]\to[0,\infty]$
with
\begin{equation}\label{eqINTFmA}
\int\limits_{\delta}^{\infty}\log\,\Phi(t)\,\frac{dt}{t^2}\ =\
+\infty\end{equation} for some $\delta>0$. Then the potential
equation (\ref{eqPotential}) has $A-$harmonic solutions of the
Dirichlet problem (\ref{eqDIR}) with the representation
(\ref{eqHYDROmA}) for each continuous function $\varphi:\partial
D\to{\Bbb R}$.}

\medskip

{\bf Corollary 24.} {\it Let $D$ be a bounded domain in ${\Bbb C}$
with no boundary component degenerated to a single point,
$A:D\to\mathbb S^{2\times 2}$ be a measurable function in $D$ such
that, for some $\alpha>0$,
\begin{equation}\label{eqEXPAmA}\int\limits_{D}e^{\alpha K_{\mu_A}(z)}\,dm(z)\ <\
\infty\ .
\end{equation} Then the potential equation
(\ref{eqPotential}) has $A-$harmonic solutions of the Dirichlet
problem (\ref{eqDIR}) with the representation (\ref{eqHYDROmA}) for
each continuous function $\varphi:\partial D\to{\Bbb R}$.}

\medskip

Thus, we have a number of effective criteria for solvability of the
Dirichlet problem to the main equation (\ref{eqPotential}) of the
hydromechanics (fluid mechanics) in strongly anisotropic and
inhomogeneous media.

\medskip

{\bf Remark 18.} By the Stoilow theorem, see e.g. \cite{Sto}, a
multi-valued solution $f=u+iv$ of the Dirichlet problem
(\ref{eqDIR}) for the Beltrami equation (\ref{eqBeltrami}) with
$K_{\mu_A}\in L^1_{\rm loc}(D)$ can be represented in the form
$f={\cal A}\circ F$ where $\cal A$ is a multi-valued analytic
function and $F$ is a homeomorphic regular solution of
(\ref{eqBeltrami}) with $\mu :=\mu_A$ in the class $W_{\rm
loc}^{1,1}$. Thus, by Theorem 5.1 in \cite{RSY$_5$}, see also
Theo\-rem 16.1.6 in \cite{AIM}, the condition (\ref{eqINTFmA}) is
not only sufficient but also necessary to have $A-$harmonic
solutions $u$ of the Dirichlet problem (\ref{eqDIR}) to potential
equations (\ref{eqPotential}) with the integral constraints
(\ref{eqINTKmA}) for all continuous functions $\varphi:\partial
D\to\Bbb{R}$, see also Remark 3.

%\medskip \medskip
{\bf \noindent Vladimir Gutlyanskii} \\
Institute of Applied Mathematics and Mechanics\\
of National Academy of Sciences of Ukraine, \\
Slavyansk 84 100,  UKRAINE\\
vgutlyanskii@gmail.com

%\medskip
%\medskip
{\bf \noindent Vladimir Ryazanov} \\
Institute of Applied Mathematics and Mechanics\\
of National Academy of Sciences of Ukraine, \\
Slavyansk 84 100,  UKRAINE\\
vl.ryazanov1@gmail.com

\medskip
%\medskip
{\bf \noindent Evgeny Sevost'yanov} \\
{\bf 1.} Zhytomyr Ivan Franko State University,  \\
Zhytomyr 10 008, UKRAINE \\
{\bf 2.} Institute of Applied Mathematics and Mechanics\\
of National Academy of Sciences of Ukraine, \\
Slavyansk 84 100,  UKRAINE\\
esevostyanov2009@gmail.com

%\medskip
%\medskip
{\bf \noindent Eduard Yakubov} \\
H.I.T. - Holon Institute of Technology,\\
Holon, PO Box 305,  ISRAEL\\
yakubov@hit.ac.il

\end{document}